\documentclass[a4paper, english, 11pt]{article}

\usepackage{setspace}

\usepackage[T1]{fontenc}
\usepackage{hyperref}
\usepackage[utf8]{inputenc}
\usepackage{dcolumn}
\usepackage{color}
\usepackage[x11names]{xcolor}
\usepackage[thmmarks]{ntheorem}
\usepackage{ragged2e}
\usepackage{pdfpages}
\usepackage{amsmath}
\usepackage[french, english]{babel}
\usepackage{multirow}
\usepackage{array}
\usepackage{booktabs}
\usepackage{wrapfig}
\usepackage{multicol}
\usepackage{algorithm}
\usepackage{lscape}
\usepackage{pifont}
\usepackage{graphicx}
\usepackage{cite}
\usepackage{calc}
\usepackage{qtree} 
\usepackage[all]{xy}
\usepackage{framed}
\usepackage{tikz}
\usepackage{mathrsfs}
\usepackage{rotating}

\def\dar[#1#2]{\ar@<2pt>[#1]\ar@<-2pt>[#2]}

\usepackage{amsfonts,amssymb,amsmath}
\usepackage[top=3cm,bottom=3cm,left=3cm,right=3cm]{geometry}
\usepackage{fancyhdr}
\usepackage{array}

\newcommand{\K}[0]{\text{K}}
\newcommand{\kk}[0]{\textbf{k}}
\newcommand{\Mod}[0]{\text{Mod}}

\newcommand{\Pe}{\text{Pers}}
\newcommand{\0}{\vec{0}}
\newcommand{\R}[0]{\mathbb{R}}
\newcommand{\N}[0]{\mathbb{N}}
\newcommand{\Q}[0]{\mathbb{Q}}
\newcommand{\Z}[0]{\mathbb{Z}}

\newcommand{\B}[0]{\mathbb{B}}
\newcommand{\F}[0]{\mathcal{F}}
\newcommand{\G}[0]{\text{R}\Gamma}

\newcommand{\V}[0]{\mathbb{V}}
\newcommand{\D}[0]{\text{D}}
\newcommand{\Obj}[0]{\text{Obj}}
\newcommand{\Ho}[0]{\text{H}}
\newcommand{\Hom}[0]{\text{Hom}}
\newcommand{\Rr}[0]{\text{R}}
\newcounter{nfigure} 
\newcommand{\cc}[0]{\mathscr{C}}
\makeatletter
 
\newskip\@bigflushglue \@bigflushglue = -100pt plus 1fil

\def\bigcentering{\let\\\@centercr\rightskip\@bigflushglue%
\leftskip\@bigflushglue
\parindent\z@\parfillskip\z@skip}

\theoremheaderfont{\scshape\bfseries}
\theorembodyfont{\normalfont}
\newtheorem{prop}{Proposition}[section] 

{ \theoremsymbol{\flushright{$\square$}} 
 \newtheorem*{pf}{Proof} }
 
 { \theoremnumbering{Roman} \theorembodyfont{\normalfont}  }
 
 \newtheorem{cor}[prop]{Corollary}
 \newtheorem{remark}[prop]{Remark}
 \newtheorem{lem}[prop]{Lemma}
 \newtheorem{defi}[prop]{Definition}
 \newtheorem{thm}[prop]{Theorem}
\author{Nicolas Berkouk, Grégory Ginot}

\begin{document}
\selectlanguage{english}

\title{A derived isometry theorem for sheaves}
\maketitle
\begin{abstract}

	Persistent homology has been recently studied with the tools of sheaf theory in the derived setting by Kashiwara and Schapira \cite{Kash18} after J. Curry has made the first link between persistent homology and sheaves.
    
    We prove the isometry theorem in this derived setting, thus expressing the convolution distance of sheaves as a matching distance between combinatorial objects associated to them that we call graded barcodes. This allows to consider sheaf-theoretical constructions as combinatorial, stable topological descriptors of data, and generalizes the situation of persistence with one parameter. To achieve so, we explicitly compute all morphisms in $\D^b_{\R c}(\kk_\R)$, which enables us to compute distances between indecomposable objects. Then we adapt Bjerkevik's stability proof to this derived setting.  
   
   As a byproduct of our isometry theorem, we prove that the convolution distance is closed, give a precise description of connected components of  $\D^b_{\R c}(\kk_\R)$ and provide some explicit examples of computation of the convolution distance.
\end{abstract}

\tableofcontents

\section{Introduction}

Persistence theory appeared in the early 2000's as an attempt to make some constructions inspired by Morse theory  computable in practice. For instance, in the context of studying the underlying topology of a data set. It has since been widely developed and applied in many ways. We refer the reader to \cite{Oudo15,Edel10} for extended expositions of the theory and of its applications. One promising expansion of the theory, initiated by Curry in his Ph.D. thesis \cite{Curr14}, is to combine the powerful theory of sheaves with ideas coming from persistence, which are driven by applications in machine learning. However, sheaf theory takes its full strength in the derived setting and Kashiwara and Schapira developed persistent homology in this new framework in \cite{Kash18}. In this paper, we show that the main theorems of one-parameter persistence theory (which we recall below) admit an analogue  in the context of derived sheaves on the real line, equipped with the convolution distance.

To our knowledge, this is the first result allowing to actually use sheaves as topological descriptors for noisy data sets.

\subsection*{One-parameter persistence}
The theory of one parameter persistence mainly relies on one construction and two theorems that we now explain. Given a real-valued function $f$ on a topological space $X$ and $i \in \Z$, consider $\mathcal{S}_i(f)(s) := H_i(f^{-1}(-\infty ; s))$ the $i$-th singular homology group of the sublevel set of $f$ with coefficient in the fixed field $\kk$. Then for $s\leq t$, the inclusion $f^{-1}(-\infty ; s)\subset f^{-1}(-\infty ; t)$ induces a linear map $\mathcal{S}_i(f)(s) \to \mathcal{S}_i(f)(t)$, and the functorial nature of singular homology gives $\mathcal{S}_i(f)$ the structure of a functor from the poset category $(\R,\leq)$ to the category  $\Mod(\kk)$ of $\kk$-vector spaces. This functor, that we still write as $\mathcal{S}_i(f)$, is usually referred to as the $i$-th sublevel-sets persistence module associated to $f$. More generally, the category $\Pe(\kk^\R)$ of \emph{persistence modules over $\R$} is precisely  the category  of functors  $(\R,\leq)\to \Mod(\kk) $.

 In \cite{Craw12}, Crawley-Boevey proved that under some finiteness assumptions on $f$, satisfied for instance when $\mathcal{S}_i(f)(s)$ is finite-dimensional at every $s\in\R$, $\mathcal{S}_i(f)$ decomposes as a locally finite direct sum of persistence modules which are constant, with values $\kk$, and supported on a given list of intervals of $\R$. This list of intervals entirely characterizes the isomorphism class of $\mathcal{S}_i(f)$ and is called the $i$-th \emph{barcode} of $f$, written $\B_i(f)$. One fundamental property that implies that the whole theory can be handled by a computer is that $\B_i(f)$ is a complete, discrete invariant of $\mathcal{S}_i(f)$.

On the other hand, for $\B_i(f)$ to be a meaningful descriptor of real-word  --hence noisy-- datasets, it must satisfy some form of stability with respect to $f$. More precisely,  it is important to understand under which distance  $\B_i(f)$ and $\B_i(g)$ are close, for $f$ and $g$ two functions $\varepsilon$-close in $L_\infty$-norm, that is in uniform convergence distance. An answer was  first given in 2005 by D. Cohen-Steiner, H. Edelsbrunner, and J. Harer in \cite{CSEH05} and is now referred to as the stability theorem. It states that if $f,g : X \to \R$ are $\varepsilon$-close in the $L_\infty$-norm, then there exists a one-to-one pairing between the intervals of $\B_i(f)$ and $\B_i(g)$, such that the right (resp. left) endpoints of each interval within a pair are closer than $\varepsilon$, and intervals can be paired to $0$ if they have length less than $2\varepsilon$. Such a pairing is called an $\varepsilon$-matching between $\B_i(f)$ and $\B_i(g)$, and we can define the bottleneck distance between $\B_i(f)$ and $\B_i(g)$ to be the infimum of the values of $\varepsilon$ for which there exists an $\varepsilon$-matching between $\B_i(f)$ and $\B_i(g)$. The stability theorem can now then be restated as follows : the bottleneck distance between $\B_i(f)$ and $\B_i(g)$ is less or equal than the $L_\infty$-norm of $f-g$.

In 2009, Chazal, Cohen-Steiner, Glisse, Guibas, and Oudot \cite{CSGG09} expressed the stability theorem algebraically, introducing the interleaving distance between one-parameter persistence modules and proving that an $\varepsilon$-interleaving (a kind of approximate isomorphism) induces an $\varepsilon$-matching between their associated barcodes. This statement is usually referred to as either the algebraic stability theorem or the \emph{isometry theorem}, and is the cornerstone of persistence techniques in the one-parameter case together with Crawley-Boevey's theorem \cite{Craw12}.

\subsection*{Persistence and sheaves}

The need for studying persistence modules obtained from functions valued in higher-dimensional  vector spaces naturally arises from the context of data analysis, see for example \cite{Lesn15,LW15}.
However, as shown in \cite{Carl09}, the category $\Pe(\kk^{\R^n})$ of functors $(\R^n,\leq) \to \Mod(\kk)$ seems to be too general for $n\geq 2$. Indeed, it contains a full sub-category equivalent to the one of finitely generated $\Z^n$-graded $\kk[x_1,...,x_n]$-modules, which implies that there is no hope for a barcode like decomposition when $n \geq 2$. There are mainly two directions undertaken to handle this issue.

The first one, initiated by Magnus Botnan and Michael Lesnick  \cite{MB17,BoLe16}, then pursued by Bjerkevik in \cite{Bjer16}, Cochoy and Oudot in \cite{CO17}, consists in restricting the study of $\Pe(\kk^{\R^n})$ to simpler sub-categories, for example, the one of persistence modules that admit a decomposition into interval modules as in the one parameter case. In \cite{Bjer16}, Bjerkevik proves that the bottleneck distance of two interval decomposable modules is bounded by a multiple (depending on the number of parameters) of the interleaving distance, and in \cite{CO17}, Cochoy and Oudot prove that a certain kind of persistence modules over $\R^2$, namely the pointwise-finite dimensional exact bi-modules, actually have an interval decomposition. Nevertheless, it remains unclear how likely is the interval decomposable case  to appear in practice.

On the other hand, it seems natural to treat persistence modules as sheaves, which is precisely what Justin Curry initiated in his Ph. D. thesis \cite{Curr14} by expressing persistence ideas in the formalism of (co-)sheaves of vector spaces on topological spaces. In particular, he defined a distance on the category of sheaves inspired by the interleaving distance, based on convolution. He also asked whether in the case of sheaves over $\R$, this distance could be expressed as a bottleneck distance. In 2018, Kashiwara and Schapira \cite{Kash18} introduced independently a derived version of the constructions of persistence theory in the category of sheaves on real vector spaces, defining the convolution distance on its derived category, proving a stability theorem and introducing a promising notion of higher-dimensional barcodes for a large category : the $\gamma$-piecewise linear sheaves.

\subsection*{Content of the paper}

In this paper, we provide answers to the question asked by Justin Curry at the end of his thesis, in the setting of Kashiwara and Schapira. We will explain later on our choice and motivations to work in this derived setting. 

%\color{blue}
It follows from general theorems that $\D^b_{\R c}(\kk_\R)$ -- the bounded derived category of constructible sheaves on $\R$ -- is a Krull-Schmidt category, whose indecomposable objects are constant sheaves over some real interval, concentrated in one degree. Recording all of the interval appearing in this decomposition, together with their degree, lead to the notion of graded barcode. Since Kashiwara and Schapira have equipped this category with the convolution distance, which is inspired by the interleaving distance, a natural question is whether the convolution distance between two sheaves in $\D^b_{\R c}(\kk_\R)$ can be computed as a matching between graded intervals, in the same fashion as the well-known isometry theorem for one parameter persistence. We provide a positive answer to this question, proving a derived isometry theorem in this setting. One important feature of our work is that the matching distance associated to the convolution distance allows to match intervals accross degree, this being due to the fundamentally derived nature of the convolution distance.  We hope that this result will open the door to considering sheaf-theoretical constructions as tools for applications in machine learning.  In particular, one future direction of research is to elucidate the implications of this derived isometry theorem with level-sets persistence, which we have undertaken in \cite{berkouk2019levelsets}. Our result might also be of independent interest for mathematical fields where barcodes techniques have allowed to obtain new results, such as symplectic topology \cite{polterovich2021topological,Asano_2020}.

%\color{black}

The paper is structured as follows : 

\begin{enumerate}
\item  Section 2 aims at introducing the mathematical context of the paper. We give a general definition of the isometry theorem problem for metric Krull-Schmidt categories. We also explain the convolution distance and the decomposition result of constructible sheaves over $\R$ obtained by Kashiwara and Schapira in \cite{Kash18}. 

\item Section 3 is dedicated to the complete description of the morphisms in $\D^b_{\R c}(\kk_\R)$, the derived category of constructible sheaves on $\R$, and to compute the action of the convolution functor $-\star \K_\varepsilon$. Note that the computations of this section (propositions \ref{P:ComputationDerivedConvolution} and \ref{P:TableHomDerived}),  may be of independent interest. 

\item Section 4 describes the conditions for two indecomposables sheaves to be $\varepsilon$-close, and introduces the notion of CLR decomposition for any sheaf $F\in \D^b_{\R c}(\kk_\R)$, which is adapted to the convolution distance in the following sense : two sheaves are $\varepsilon$-close with respect to $d_C$ if and only if their central (resp. left, resp. right) parts are. 

\item In Section 5, we prove that given an $\varepsilon$-interleaving between the central (resp. left, resp. right) parts of two sheaves, it induces an $\varepsilon$-matching between the graded-barcodes of their central (resp. left, resp. right) parts. We reduce the proof for left and right parts to the well-known case of one-parameter persistence modules by introducing a family of fully faithful functors from sheaves supported on half-open intervals to persistence modules. The construction of the $\varepsilon$-matching between the central parts is far less direct. We adapt the proof of Bjerkevik \cite{Bjer16} to our setting, introducing a similar pre-order $\leq_\alpha$ on central parts, enabling us to \og trigonalize \fg ~ the interleaving morphisms. By a rank argument, this allows us to apply Hall's marriage theorem and to deduce the existence of a $\varepsilon$-matching. Note that our definition of $\leq_\alpha$ differs in nature from Bjerkevik's, for it enables us to compare elements of the graded-barcodes in \emph{different degrees}. We conclude the section by proving the isometry theorem~\ref{T:DerivedIsometry}, which states that \og $ d_C =  d_B $ \fg.

\item Section 6 provides some applications of the isometry theorem. We start by an example with explicit computations brought to our knowledge by Justin Curry and that motivated our work. Then, we prove that the convolution distance is closed (two sheaves are $\varepsilon$-close if and only if they are $\varepsilon$-interleaved) thus answering an open question of \cite{Kash18} in dimension one. We provide a counter-example of two non constructible sheaves $F$ and $G$ such that $d_C(F,G) = 0$ but $F \not \simeq G$.  We also prove that the open balls of the metric space $(\D^b_{\R c}(\kk_\R),d_C)$ are path-connected, hence leading to a characterization of connected components of $\D^b_{\R c}(\kk_\R)$ (Theorem~\ref{T:pathcomponent}). To do so, we introduce an explicit skeleton for  $\D^b_{\R c}(\kk_\R)$ called the category  $Barcode$ of barcodes (Definition~\ref{D:CatBarcode}). The latter inherits an extended metric space structure.
\end{enumerate}

\paragraph{Acknowledgement} The authors would like to thank Justin Curry, Steve Oudot, Pierre Schapira and Magnus Botnan for many enlightening discussions with them. 

The first author is supported by Innosuisse grant 45665.1 IP-ICT, the second author was partially supported 
by ANR grants Catag and Chrok.

\section{Preliminaries}\label{S:Preliminaries}
This section aims at introducing the notation that we will use in this paper, presenting the theoretical framework of \cite{Kash18} and explaining precisely the problem underlying the isometry theorem. 

\subsection{The isometry theorem problem for metric Krull-Schmidt categories} \label{S:IsomProblem}

%\color{blue}
Let $\cc$ be an additive category. Recall that an object $M \in \cc$ is indecomposable if $M \not = 0$, and for any isomorphim $M \simeq M_1 \oplus M_2$, either $M_1$ or $M_2$ is equal to $0$.

\begin{defi} \label{D:KrullSchmidt}
A category $\cc$ is a {Krull-Schmidt category} if it satisfies the following axioms. 
\begin{itemize}
    \item[(KS-1)] $\cc$ is an additive category.
    \item[(KS-2)] For any object $X$ of $\cc$, there exists a family of indecomposable objects $\mathbb{B}(X)$ of $\cc$ such that  $X \simeq \bigoplus_{I \in \mathbb{B}(X)} I  $ which is essentially unique. That is, for any other family of indecomposable objects $\mathbb{B}'(X)$ with the same property, there exists a bijection $\sigma:\mathbb{B}(X) \to \mathbb{B}'(X)$ such that $I\simeq \sigma(I)$, for all $I$ in $\mathbb{B}(X)$.
     \item[(KS-3)] For any object $X$ of $\cc$ such that  $ X \simeq \bigoplus_{I \in \mathbb{B}(X)} I  $ with $\mathbb{B}(X)$ a collection of indecomposable objects of $\cc$, then $\prod_{I \in \mathbb{B}(X)} I $ exists in $\cc$ and the canonical morphism : $$\bigoplus_{I \in \mathbb{B}(X)} I \longrightarrow \prod_{I \in \mathbb{B}(X)} I  $$ is an isomorphism. 
    
\end{itemize}
    
\end{defi}

\begin{remark}
{Note that the usual definition of a Krull-Schmidt category asks, with notations of definition \ref{D:KrullSchmidt}, that $\mathbb{B}(X)$ is finite (see \emph{eg.} \cite{Krau14}). This will not be sufficient for our study of constructible sheaves over $\R$ since they are potentially infinite direct sum of sheaves constant on a real interval (theorem \ref{T:KSdecomposition}). However, one important behaviour of these direct sums is that they satisfy axiom (KS-3), since they are locally finite. }
\end{remark}

\begin{defi}\label{D:metricCategory}
Let $\cc$ be any category. An extended pseudo-distance on $\cc$ is a map $d$ defined on  $\Obj(\cc) \times \Obj(\cc)$ satisfying, for all $X,Y,Z$ in $\cc$: 
\begin{enumerate}
\item[(M1)] $d(X,Y)\in \R_{\geq 0}\cup \{+\infty\}$,
\item[(M2)] $d(X,Y) = d(Y,X)$,
\item[(M3)] $d(X,Z) \leq d(X,Y) + d(Y,Z)$,
\item[(M4)] if $X \simeq Y$, then $d(X,Y)=0$.
\end{enumerate}
In this situation, $(\cc,d)$ will be called a \emph{metric category}.
\end{defi}

Let $(\cc,d)$ be a metric category, such that $\cc$ is a Krull-Schmidt category. For any object $X$ of $\cc$, one denotes by $\mathbb{B}(X)$ a collection of indecomposables objects of $\cc$ such that $ X\simeq  \bigoplus\limits_{I\in \mathbb{B}(X)} I$. 

\begin{defi} \label{D:matching}
Let $X,Y \in \cc$ and $\varepsilon\geq 0$. An {$\varepsilon$-matching} between $\mathbb{B}(X)$ and $\mathbb{B}(Y)$ is the following data: two subcollections $\mathcal{X} \subset \mathbb{B}(X)$ and $\mathcal{Y}\subset \mathbb{B}(Y)$, and a bijection $\sigma:\mathcal{X} \to \mathcal{Y} $ satisfying: 
\begin{enumerate}
    \item  $d(I,\sigma(I)) \leq \varepsilon$ for all $I$ in $\mathcal{X}$;
    \item  $d(I, 0) \leq \varepsilon$, for all $I$ in $\mathbb{B}(X) \backslash \mathcal{X}$ or in $\mathbb{B}(Y) \backslash \mathcal{Y}$.
\end{enumerate}
\end{defi}

In this situation, we will use the notation $\sigma : X \not \to Y$ and designate $\mathcal{X}$ by $\text{coim}(\sigma)$ (resp. $\mathcal{Y}$ by $\text{im}(\sigma)$).

Since $d$ satisfies $(M4)$, the existence of a $\varepsilon$-matching does not depend on the choice of representatives in $\mathbb{B}(X)$ and $\mathbb{B}(Y)$. 

\begin{defi} \label{D:bottleneckgeneral} Let $X$ and $Y$ be two objects of $\cc$. One defines the bottleneck distance associated to $d$ between $X$ and $Y$ as the possibly infinite following quantity:
$$d_B(X,Y) = \inf \{\varepsilon \geq 0 \mid \text{there exists an $\varepsilon$-matching between $\mathbb{B}(X)$ and $\mathbb{B}(Y)$}\}. $$

\end{defi}

\begin{prop}
The map $d_B$ is an extended pseudo-metric on $\cc$.
\end{prop}

\begin{pf} The fact that $d_B$ takes values in $\R_{\geq 0}\cup \{+\infty\}$ and Properties  (M2) and (M4) are directly inherited from $d$. For the triangle inequality, observe that since $d$ satisfies (M3), we can compose a $\varepsilon$-matching between $\B(X)$ and $\B(Z)$, with a $\varepsilon'$-matching between $\B(Z)$ and $\B(Y)$ to obtain a $\varepsilon + \varepsilon'$ matching between $\B(X)$ and $\B(Y)$.

\end{pf}

\textbf{Terminology.} The {isometry theorem problem} associated to the metric Krull-Schmidt category $(\cc,d)$ is to determine whether $d=d_B$. 

\color{black}

\subsection{Notations for sheaves and complexes}

Throughout the paper and except when stated otherwise, we will follow the notations introduced in \cite{Kash18} and \cite{Kash90}. We will also freely refer to some of their proofs.

In the paper, $\kk$ will denote a field, $\Mod(\kk)$ the category of vector spaces over $\kk$, $\text{Mod}_f(\kk)$ the category of finite dimensional vector spaces over $\kk$. Let $X$ be a topological space. Then we will note $\Mod(\kk_X)$ the category of sheaves of $\kk$-vector spaces on $X$. For shortness, we will also write $\Hom$ for $\Hom_{\Mod(\kk_\R)}$.  

For $\cc$ an abelian category, denote $\text{C}^b(\cc)$ its category of bounded complexes, $\text{K}^b(\cc)$ its bounded homotopy category and $\D^b(\cc)$ its bounded derived category. For simplicity, we shall write $\D^b(\kk)$ instead of $\D^b(\Mod(\kk))$ and $\D^b(\kk_X)$ instead of $\D^b(\Mod(\kk_X))$. When the context is clear, we will simply call  sheaves the objects of $\D^b(\kk_X)$. For a complex $F \in \text{C}^b(\cc)$ and an integer $k$, define the $k$-th shift of $F$ by: for $n\in \Z$, $F[k]^n = F^{k+n}$ and $d_{F[k]}^n= (-1)^{k} d_F^{n+k}$.

We will use the classical notations of \cite{Kash90} for the Grothendieck operations on sheaves. Moreover, we recall the following : for $X_1$ and $X_2$ two topological spaces, we denote $p_i: X_1 \times X_2 \to X_i$, $i=1,2$ the canonical projections. Let  $F \in \Mod(\kk_{X_1})$ and $G\in  \Mod(\kk_{X_2})$, define their {external tensor product} $F \boxtimes G \in \Mod(\kk_{X_1 \times X_2})$ by the formula : $$F \boxtimes G := p_1^{-1}F \otimes p_2^{-1}G.$$

Observe that since we are working over a field, this operation is exact, hence need not to be derived.

\begin{defi}
For $M$ a real analytic manifold, and $F\in \Mod(\kk_M)$, $F$ is said to be {weakly $\R$-constructible} if there exists a locally finite sub-analytic stratification of $M = \sqcup_\alpha M_\alpha$, such that for each stratum $M_\alpha$, the restriction $F_{|M_\alpha}$ is locally constant. If in addition, the stalks $F_x$ are of finite dimension for every $x\in M$, we say that $F$ is {$\R$-constructible}. We might often say constructible instead of $\R$-constructible, since, in this paper, it is the only notion of constructibility we use.
\end{defi}

We will write $\Mod_{\R c}(\kk_M)$ for the abelian category of $\R$-constructible sheaves on $M$, and $\D^b_{\R c}(\kk_M)$ the full triangulated subcategory of $\D^b(\kk_M)$ consisting of complexes of sheaves whose cohomology objects lie in $\Mod_{\R c}(\kk_M)$. Note that Theorem 8.4.5 in \cite{Kash90} asserts that the natural functor $\D^b(\Mod_{\R c}(\kk_M)) \to \D^b_{\R c}(\kk_M)$ is an equivalence of triangulated categories.

\subsection{Constructible sheaves over $\R$}
  Theorem~\ref{T:KSdecomposition} below is proved in \cite{Kash18} and generalizes Crawley-Boeyvey's theorem \cite{Craw12} to the context of constructible sheaves on the real line. Together with Theorem~\ref{T:KSstructure}, they will be the cornerstone to prove that $\D^b_{\R c}(\kk_\R)$ is a Krull-Schmidt category. 

\begin{defi}
Let $\mathcal{I}=\{I_\alpha\}_{\alpha \in A}$ be a multi-set of intervals of $\R$, that is, a list of interval where one interval can appear several times. Then $\mathcal{I}$ is said to be {locally finite} if and only if for every compact set $K\subset \R$, the set $\{\alpha \in A \mid K  \cap I_\alpha \not = \emptyset\}$ is finite.
\end{defi}

\begin{thm}[Theorem 1.17 - \cite{Kash18}]\label{T:KSdecomposition}
Let $F \in \Mod_{\R c}(\kk_\R)$, then there exists a unique locally finite multi-set of intervals $\mathbb{B}(F)$ such that $$F \simeq \bigoplus_{I \in \mathbb{B}(F)} \kk_{I}.$$ Moreover, this decomposition is unique up to isomorphism.
\end{thm}

\begin{defi}\label{D:Barcode}
The multi-set $\B(F)$ is the called the barcode of $F$.
\end{defi}

\begin{cor}
Let $F,G \in \D^b_{\R c}(\kk_\R)$, and $j\geq 2$, then: $\text{Ext}^j(F,G) = 0.$

\end{cor} 

A classical consequence of such a statement is the following: 

\begin{thm}\label{T:KSstructure}
Let $F\in \D^b_{\R c}(\kk_\R)$. Then there exists an isomorphism in $\D^b_{\R c}(\kk_\R)$: $$F \simeq \bigoplus_{j\in\Z} \Ho^j(F)[-j]$$
where $\Ho^j(F)$ is seen as a complex concentrated in degree $0$. Not that this isomorphism is not functorial.
\end{thm}

\begin{defi}\label{D:GradedBarcodes}
Let $F\in \D^b_{\R c}(\kk_\R)$, we define its \emph{graded-barcode} $\B(F)$ as the collection $(\mathbb{B}^j(F))_{j\in \Z}$ where $\B^j(F) := \B(\Ho^j(F))$. Furthermore, to indicate that an interval $I\subset \R$ appears in degree $j\in \Z$ in the graded-barcode of $F$, we will write $I^j \in \B(F)$. The element $I^j$ is called a graded-interval.
\end{defi}

Note that  Theorems \ref{T:KSdecomposition} and \ref{T:KSstructure} imply that $$F \simeq \bigoplus_{\substack{I^j \in \B(F)}} \kk_I[-j].  $$
Therefore, $\B(F)$ is a complete discrete invariant of the isomorphism class of $F$ in $\D^b_{\R c}(\kk_\R)$.

\begin{prop}\label{P:SheafKS}
The category $\D^b_{\R c}(\kk_\R)$ is Krull-Schmidt.
\end{prop}

\begin{pf}
The category $\D^b_{\R c}(\kk_\R)$ is additive by construction. Moreover, for $I \subset \R$ an interval and $j\in \Z$, the sheaf $\kk_I[-j]$ is indecomposable in $\D^b_{\R c}(\kk_\R)$. Therefore, we deduce from Theorems \ref{T:KSdecomposition} and \ref{T:KSstructure} that $\D^b_{\R c}(\kk_\R)$ satisfies (KS2). Let us now prove that $\D^b_{\R c}(\kk_\R)$ satisfies (KS3). Let $F \in \D^b_{\R c}(\kk_\R)$,  then there exists $N\geq 0$ such that $\B^j(F) = \emptyset$ for $|j| \geq N$. As a consequence,  given  $x \in \R$, the set $\{I \in \sqcup_j \mathbb{B}^j(F) \mid (x-1,x+1) \cap I \not = \emptyset\}$ is finite. We deduce that the natural morphism: $$ F \simeq \bigoplus_{\substack{I^j \in \B(F)}} \kk_I[-j] \longrightarrow \prod_{\substack{I^j \in \B(F)}} \kk_I[-j]$$
induces an isomorphism when taking the stalk at $x$. This being true for all $x \in \R$, it is an isomorphism.
\end{pf}

\subsection{Metric for sheaves}

In \cite{Curr14}, Curry defined an interleaving-like distance on $\Mod(\kk_X)$, for $(X,d)$ a metric space.
It is based on what he calls the smoothing of opens. For $F\in \Mod(\kk_X) $, define $F^\varepsilon \in \Mod(\kk_X)$ as
the sheafification of $U \mapsto F(U^\varepsilon)$, 
with $U^\varepsilon = \{x\in X \mid \exists u \in U, d(x,u) \leq \varepsilon\}$.
This yields a functor $[\varepsilon] : \Mod(\kk_X) \to \Mod(\kk_X) $
together with a natural transformation $[\varepsilon] \Rightarrow \text{id}_{\Mod(\kk_X)}$.
Although this seems to mimic the construction of interleaving distance for persistence modules, 
one must pay attention to the fact  that $[\varepsilon]$ is only left-exact. Since topological informations are obtained 
from sheaves by considering sheaf-cohomology, one needs to derive the functor $[\varepsilon]$
in order to keep track of cohomological informations while smoothing a sheaf. 
This is precisely the sense of the construction of Kashiwara and Schapira using convolution of sheaves, 
which has the advantage to have a nice expression in term of Grothendieck operations (that allows appropriate operations  for sheaf cohomology). 

In this section, we make a short review of the concepts introduced in \cite{Kash18}. 
The framework is the study of sheaves
on a real vector space $\V$  of finite dimension $n$ equipped with a norm $\| \cdot \|$. 
For two such sheaves, one can define their convolution, which,
as the name suggests, will be at the core of the definition of the convolution distance (definition \ref{D:convolutionDistance}).

The construction of the convolution of sheaves is as follows.
Consider the following maps (addition and the canonical projections): 

$$s : \V \times \V \to \V, ~~~s(x,y) = x + y, $$ 
$$q_i : \V \times \V \to \V ~~(i=1,2) ~~~q_1(x,y) = x,~q_2(x,y) = y.$$

\begin{defi}\label{D:Convolution}
For $F,G\in \D^b(\kk_\V)$, we define the {convolution} of $F$ and $G$ by the formula: $$F\star G = \text{R}s_!(F\boxtimes G).$$

\end{defi}

This defines a bi-functor : $\D^b(\kk_\V)\times \D^b(\kk_\V) \to \D^b(\kk_\V)$. 

In the following, we will be interested in a more specific case : the convolution will be considered with one of the sheaves being the constant sheaf supported on a ball centered at 0. 

For $r >0$, we denote $B_r := \{x \in \V \mid \|x\| \leq r\}$ the closed ball of radius $r$ centered at $0$, and $\overset{\circ}{B}_r$ its interior, that is, the open ball of radius $r$ centered at 0. For $\varepsilon \in \R$ we define : 

\begin{equation}\label{eq:ballkernel}
    \K_\varepsilon := \begin{cases}
    \kk_{B_\varepsilon} ~\text{if}~\varepsilon \geq 0 \\
    \kk_{\overset{\circ}{B}_{-\varepsilon}}[\dim(\V)] ~\text{if}~\varepsilon<0
    \end{cases} 
\end{equation} 
with $\kk_{\overset{\circ}{B}_{-\varepsilon}}[\dim(\V)]$, seen as a complex concentrated in degree $-\dim(\V)$. We have the following properties: 

\begin{prop}[Section 2.1 - \cite{Kash18}]\label{P:propertiesofconvolution} Let $\varepsilon, \varepsilon'\in \R$ and $F \in \D^b(\kk_\V)$ .
\begin{enumerate}

\item One has functorial isomorphisms $(F\star K_{\varepsilon} )\star K_{\varepsilon'} \simeq F \star K_{\varepsilon + \varepsilon'} $ and $F\star K_0\simeq F$.

\item If $\varepsilon' \geq \varepsilon $, there is a canonical morphism 
$K_{\varepsilon'}\to K_{\varepsilon}$ in $\D^b(\kk_\V)$
inducing a natural transformation $F\star K_{\varepsilon'} \to F \star K_{\varepsilon} $. 
In the special case where $\varepsilon = 0$, we shall write $\phi_{F, \varepsilon'}$ for this natural transformation.

\item The canonical morphism $F\star K_{\varepsilon'} \to F \star K_{\varepsilon}$ induces an isomorphism
\[\G(\V; F\star K_{\varepsilon'}) \tilde{\to} \G(\V; F\star \K_\varepsilon).\]  
\end{enumerate}

\end{prop}

\begin{defi}\label{D:smoothingmorphism}
With the same notations as in the previous proposition,  the morphism $\phi_{F,  \varepsilon'}$ is called the $\varepsilon'$-smoothing morphism of $F$.

\end{defi}

In particular, Proposition~\ref{P:propertiesofconvolution} implies that any map $f: F\star K_{\varepsilon} \to G$ 
induces canonical maps  
\begin{equation} \label{eq:propertiesofconvolution} f\star K_{\tau}: F\star K_{\varepsilon+\tau}\simeq F\star  K_{\varepsilon}\star K_{\tau}  \to G \star K_{\tau}.\end{equation} 

The following definition is central.
\begin{defi}[Definition 2.2 -  {\cite{Kash18}}]\label{D:inteleavingconvolution}
For $F,G\in \D^b(\kk_\V)$ and $\varepsilon\geq 0$, one says that $F$ and $G$
are $\varepsilon$-{interleaved} if there exists two morphisms  $f : F\star \K_\varepsilon \to G$ 
and $g :  G \star \K_\varepsilon \to F$ (in $\D^b(\kk_\V)$) such that the compositions
$F\star K_{2\varepsilon} \stackrel{ f\star \K_\varepsilon}{\longrightarrow}\K_\varepsilon\star G \stackrel{g}{\longrightarrow} F $ and $G\star K_{2\varepsilon} \stackrel{g \star \K_\varepsilon}{\longrightarrow}\K_\varepsilon\star F \stackrel{f}{\longrightarrow} G $ are the natural morphisms $F\star K_{2\varepsilon} \stackrel{\phi_{F,2\varepsilon}}{\longrightarrow} F$ and $G\star K_{2\varepsilon} \stackrel{\phi_{G,2\varepsilon}}{\longrightarrow} G$, that is, we have a commutative diagram in $\D^b(\kk_\V)$ :

$$ \xymatrix{
F\star K_{2\varepsilon}  \ar[rd]\ar@/^0.7cm/[rr]^{\phi_{F,2\varepsilon}} \ar[r]^{f\star \K_\varepsilon} & G\star \K_\varepsilon  \ar[rd] \ar[r]^{g} & F\\
G \star K_{2\varepsilon}  \ar[ur]\ar@/_0.7cm/[rr]_{\phi_{G,2\varepsilon}} \ar[r]^{g\star \K_\varepsilon} & F\star \K_\varepsilon  \ar[ur] \ar[r]^{f} & G
   } $$ 
In this case, we write $F \sim_\varepsilon G$.
\end{defi}

Observe that $F$ and $G$ are $0$-interleaved if and only if $F\simeq G$.

\begin{remark}
One must be aware that in \cite{Kash18}, the authors call this data an $\varepsilon$-isomorphism.
Here, we choose to follow the usual terminology of persistence theory.
\end{remark}

Since 0-interleavings are isomorphisms, the existence of an $\varepsilon$-interleaving between two sheaves expresses a notion of closeness. This leads the authors of \cite{Kash18} to define the convolution distance as follows:

\begin{defi}[Definition 2.2 - \cite{Kash18}]\label{D:convolutionDistance}
For $F,G\in \D^b(\kk_\V)$, we define their {convolution distance} as: $$d_C(F,G) : = \text{inf} \left (\{ +\infty \} \cup \{a \in \R_{\geq 0 } \mid \text{$F$ and $G$ are $a$-interleaved} \} \right )$$
\end{defi}

\begin{prop}[Section 2.2 - \cite{Kash18}] \label{P:propertiesconvdistance}
The convolution distance is an extended pseudo-distance on $\D^b(\kk_\V)$ that is, it satisfies (M1)-(M4) from definition \ref{D:metricCategory}.
\end{prop}

The following proposition expresses that the functors $\G(\V;-)$ and $\G_c(\V;-)$ 
define some necessary conditions for two sheaves to be at finite convolution distance. 
This is similar to the case of interleaving distance for persistence modules $M :(\R^n,\leq) \to \Mod(\kk) $, 
where the role of $\G(\V;-)$ is played by the colimit functor over $\R^n$.
\begin{prop}[Remark 2.5 - \cite{Kash18}]
Let $F,G \in \D^b(\kk_\V)$.
\begin{enumerate}
\item If $d_C(F,G) < +\infty $ then $\G(\V,G) \simeq \G(\V,F)$ and $\G_c(\V;G) \simeq \G_c(\V;F)$.

\item If supp$(F)$, supp$(G) \subset B_a$ then $d_C(F,G) \leq 2a $ if and only if $\G(\V,G) \simeq \G(\V;F)$.
\end{enumerate}
\end{prop}

There is a fundamental example to keep in mind in the context of sheaves. This example is the one mimicking 
the persistence modules $\mathcal{S}_i(f)$ of a continuous map:   
given $X$ a topological space and $u : X \to \V$ a continuous map,
one can consider the sheaves $\Rr u_* \kk_X$ and $\Rr u_! \kk_X$. 
Roughly speaking and under some smoothness assumptions on $X$ and $f$, they contain the information on how the cohomologies 
of the fibers of $u$ evolve when moving on $\V$. 
For this information to be meaningful for applications in machine learning, it has to be stable when we perturb $u$, that is, $Ru_* \kk_X$
must stay in a neighborhood in the sense of the convolution distance, controlled by the size of the perturbation of $u$. 
This is what expresses the following theorem, which is the analogous of the stability theorem in the context of persistence theory.

\begin{thm}[Theorem 2.7 - \cite{Kash18}]\label{T:DerivedStability}
Let $X$ a locally compact topological set, and $u,v : X \to \V$ two continuous functions. Then for any $F\in \D^b(\kk_X)$ one has : $$d_C(\Rr u_*F, \Rr v_*F) \leq \| u - v\| ~~~ \text{and}~~~ d_C(\Rr u_{!}F, \Rr v_{!}F) \leq \| u - v\| $$
where we define $\|u-v \| = \sup_{x\in X} \|u(x) - v(x) \|$.
\end{thm}

\subsection{The isometry theorem problem for $(\D^b_{\R c}(\kk_\R), d_C)$}

In the previous sections, we have shown that $(\D^b_{\R c}(\kk_\R), d_C)$ is a metric Krull-Schmidt category, see Propositions \ref{P:SheafKS} and  \ref{P:propertiesconvdistance}. 
We will prove later (see~\ref{T:DerivedIsometry}) that the isometry theorem holds in this context, that is, the convolution distance equals its associated bottleneck distance. In other words, we can compute the convolution distance between two sheaves $F,G \in \D^b_{\R c}(\kk_\R)$ as a matching distance between the multiset of graded-intervals appearing in their decomposition, which we will later call their graded-barcodes (see Definition~\ref{D:Barcode}). This matching between barcodes will be similar to the usual matching between barcodes (that is in terms of comparing end points of bars and overall length) but will take into account shift of degrees and types of bars as well.

\color{black}

\section{Computations in $\D^b_{\R c}(\kk_\R)$}

This section aims at making explicit all the computations of morphisms in $\D^b_{\R c}(\kk_\R)$, and determine the action of the functor $- \star \K_\varepsilon$. 
Combining Theorems \ref{T:KSdecomposition} and \ref{T:KSstructure}, we see that any object of $\D^b_{\R c}(\kk_\R)$ is isomorphic to a direct sum of sheaves constant on an interval seen as a complex concentrated in one degree. Hence, to give a full description of the morphisms, it is enough to compute  $\Rr\Hom_{\Mod(\kk_\R)}(\kk_I,\kk_J)$ for $I,J$ two intervals -- in the sequel, we shall write for short $\Hom$ for $\Hom_{\Mod(\kk_\R)}$. Indeed, for $F,G \in \D^b_{\R c}(\kk_\R)$ we have by local finiteness of the barcodes of $F$ and $G$: 

$$\Rr \Hom (F,G) \simeq   \bigoplus_{\substack{ I^i \in \B(F) \\ J^j \in \B(G) }} \Rr \Hom(\kk_I,\kk_J)[i-j],$$

and we recall the classical formula $\Hom_{\D^b(\kk_\R)}(F,G) \simeq \Rr ^0 \Hom (F,G)$.

Our approach\footnote{which was suggested by the referee, that we wish to thank.} will rely on the use of the duality functor $\D$ introduced, for instance, in \cite[3.1.16]{Kash90}, that we will quickly review. Regarding the computations of the convolution with $\K_\varepsilon$, we will make use of a classical lemma for convolution of sheaves, that we will also review.

\subsection{Duality and morphisms in $\D^b_{\R c}(\kk_\R)$ }

Since $\R$ is a smooth manifold, it has finite c-soft dimension, and the dualizing complex $\omega_\R := a^! \kk_{\{pt\}} \simeq \kk_\R [1]$ is well defined, where $a : \R \to \{pt\}$. One defines the triangulated functor 
$ \D : \D^b(\kk_\R)^\text{op} \to \D^b(\kk_\R)  $, by $\D(F) = \Rr \mathscr{H}om(F, \omega_\R)$.
We recall the following proposition, which is a weaker version of \cite[3.4.3 and 3.4.6]{Kash90} adapted to our setting, since constructible sheaves are in particular cohomologically constructible. 

\begin{prop} \label{P:duality}
Let $F,G\in \D^b_{\R c}(\kk_\R)$, then: \begin{enumerate}
    \item the canonical map $F \stackrel{\sim}\to \D \circ \D (F)$ is an isomorphism,
    \item $\Rr \mathscr{H}om(F, G) \simeq \D(\D(G) \otimes F).  $ 
\end{enumerate} 
\end{prop}

For an open subset $U \subset \R$, one has $\D(\kk_U) \simeq \kk_{\overline{U}}[1]$, where $\overline{U}$ is the closure of $U$,  from which we deduce the following computations: 

\begin{prop}
Let $a<b \in \R$, then: 
\begin{enumerate}
    \item $\D(\kk_{\{a\}}) \simeq \kk_{\{a\}}$,
    \item $\D(\kk_{[a,b)}) \simeq \kk_{(a,b]}[1].$
\end{enumerate}
\end{prop}

\begin{pf}
To prove these results, one simply applies successively the contravariant triangulated functor $\D$ to the distinguished triangles \[\kk_{\R\backslash \{a\}} \to \kk_\R \to \kk_{\{a\}} \stackrel{+1}{\to} \quad \mbox{  and }\quad \kk_{(a,b)} \to \kk_{(a,b]} \to \kk_{\{a\}} \stackrel{+1}{\to}. \]. %Note that $\D$ is contravariant, therefore it reverses the graduation.
\end{pf}

By definition, given $F,G\in \D^b_{\R c}(\kk_\R)$, one has $\Rr \Hom (F,G)\simeq \Rr \Gamma (\R ; \Rr \mathscr{H}om(F, G))$. We now review the computations of derived global sections of indecomposable constructible sheaves on $\R$.

\begin{prop}\label{P:sections} Let $a<b \in \R$, then: 

\begin{enumerate}
    \item $\Rr \Gamma (\R ; \kk_{[a,+\infty)}) \simeq \Rr \Gamma (\R ; \kk_{(-\infty,b]}) \simeq \Rr \Gamma (\R ; \kk_{[a,b]}) \simeq \Rr \Gamma (\R ; \kk_{\{a\}}) \simeq \Rr \Gamma (\R ; \kk_\R) \simeq \kk,$
    
    \item $\Rr \Gamma (\R ; \kk_{(a,+\infty)}) \simeq \Rr \Gamma (\R ; \kk_{(-\infty,b)}) \simeq  0,$

    \item $\Rr \Gamma (\R ; \kk_{(a,b)}) \simeq \kk[-1],$
    
      \item $\Rr \Gamma (\R ; \kk_{[a,b)}) \simeq \Rr \Gamma (\R ; \kk_{(a,b]}) \simeq 0.$
\end{enumerate}

\end{prop}

\begin{pf}
\begin{enumerate}
    \item Let $I$ be any of the interval $[a,+\infty), (-\infty,b], [a,b],\{a\} $ and $i : I \to \R$ be the inclusion. Then by definition $\kk_I = i_\ast i^{-1}\kk_\R$, and we have the following isomorphisms: 
    
    \begin{align*}
        \Rr \Gamma (\R ; \kk_I) & \simeq \Rr \Hom (\kk_\R ; i_\ast i^{-1} \kk_\R) \\
        & \simeq \Rr \Hom (i^{-1}\kk_\R ;  i^{-1} \kk_\R) \\
        & \simeq \Rr \Gamma (I;   \kk_I) \\
        & \simeq \kk.
    \end{align*}

    \item This is a consequence of 1. by applying  the triangulated functor $\Rr \Gamma (\R ; - )$ to the distinguished triangles $\kk_{(a,+\infty)} \to \kk_\R \to \kk_{(-\infty,a]} \stackrel{+1}{\to}$ and $\kk_{(-\infty,b)} \to \kk_\R \to \kk_{[b,+\infty)} \stackrel{+1}{\to}$.
    
    \item This is a consequence of 2. by applying the triangulated functor $\Rr \Gamma (\R ; - )$ to the distinguished triangle $\kk_{(a,b)} \to \kk_\R \to \kk_{(-\infty,a]} \oplus \kk_{[b,+\infty)} \stackrel{+1}{\to}.$
    
    \item This is a consequence of 3. by applying the triangulated functor $\Rr \Gamma (\R ; - )$ to the distinguished triangles $\kk_{[a,b)} \to \kk_{[a,b]} \to \kk_{\{b\}} \stackrel{+1}{\to}$ and $\kk_{(a,b]} \to \kk_{[a,b]} \to \kk_{\{a\}} \stackrel{+1}{\to}$.
\end{enumerate}
\end{pf}

\begin{remark} Let $a \leq b \in  \R$. For simplicity, we adopt the following convention. When stating results about the interval $(a,b)$, we will always implicitly assume that $a<b$ .

\end{remark}

\begin{prop}\label{P:TableHomNonDerived}\label{P:TableHomDerived}
Let $a\leq b$ and $c\leq d$ in $\R$. We have the following derived morphism groups,
where the first column defines the support of the  left-side object (i.e. the source) in $\Rr \Hom(-,-)$ and the first line the right-side one : 

\newpage

\thispagestyle{empty}
\tiny
%\newgeometry{top=0.5cm}

\begin{center}
\begin{sideways}

\begin{math}

\begin{array}{l|lllllllll}

\Rr \Hom(\kk_I , \kk_J ) & (c,d)&  [c,d]  &  [c,d) &  (c,d] & (-\infty, d) & (c, +\infty) & (-\infty, d] & [c, +\infty) & \R \\ 

\hline

 (a,b) & \begin{cases} \kk ~\text{if}~ (a,b)\subset (c,d) \\ \kk [-1] ~\text{if}~ [c,d ]\subset (a,b) \\ 0~\text{else} \end{cases}  & \begin{cases} \kk ~\text{if}~ (a,b)\cap (c,d) \not = \emptyset  \\ 0~\text{else} \end{cases}   &  \begin{cases} \kk ~\text{if}~ c < b \leq d   \\ 0~\text{else} \end{cases} & \begin{cases} \kk ~\text{if}~ c \leq a < d   \\ 0~\text{else} \end{cases} & \begin{cases} \kk ~\text{if}~ b \leq d   \\ 0~\text{else} \end{cases} & \begin{cases} \kk ~\text{if}~ c \leq a   \\ 0~\text{else} \end{cases} & \begin{cases} \kk ~\text{if}~ a < d   \\ 0~\text{else} \end{cases} & \begin{cases} \kk ~\text{if}~ c < b   \\ 0~\text{else} \end{cases} & \kk  \\

 [a,b] & \begin{cases} \kk[-1] ~\text{if}~ [a,b] \cap [c,d] \not = \emptyset  \\ 0~\text{else} \end{cases} & \begin{cases} \kk ~\text{if}~ [c,d] \subset [a,b]   \\ \kk[-1] ~\text{if}~ [a,b] \subset (c,d)  \\ 0~\text{else} \end{cases}  & \begin{cases} \kk[-1] ~\text{if}~ c < a \leq d   \\ 0~\text{else} \end{cases} & \begin{cases} \kk[-1] ~\text{if}~ c \leq b < d   \\ 0~\text{else} \end{cases}  & \begin{cases} \kk[-1] ~\text{if}~ a \leq  d   \\ 0~\text{else} \end{cases}  & \begin{cases} \kk[-1] ~\text{if}~ c \leq  b   \\ 0~\text{else} \end{cases}  & \begin{cases} \kk[-1] ~\text{if}~ b < d   \\ 0~\text{else} \end{cases} & \begin{cases} \kk[-1] ~\text{if}~ c < a   \\ 0~\text{else} \end{cases} &  \kk[-1] \\

 [a,b) & \begin{cases} \kk[-1] ~\text{if}~ a \leq d <b    \\ 0~\text{else} \end{cases}   &  \begin{cases} \kk ~\text{if}~ a \leq c < b    \\ 0~\text{else} \end{cases}  & \begin{cases} \kk ~\text{if}~ a \leq c \leq b \leq d    \\ 0~\text{else} \end{cases}   & 0 & \begin{cases} \kk[-1] ~\text{if}~ a \leq d <b    \\ 0~\text{else} \end{cases}  & 0 &  0 & \begin{cases} \kk ~\text{if}~ a \leq c < b    \\ 0~\text{else} \end{cases} & 0\\

(a,b] & \begin{cases} \kk[-1] ~\text{if}~ a < c \leq b    \\ 0~\text{else} \end{cases} & \begin{cases} \kk ~\text{if}~ a < d \leq b    \\ 0~\text{else} \end{cases}  & 0 & \begin{cases} \kk ~\text{if}~ c \leq a \leq d \leq b    \\ 0~\text{else} \end{cases}  & 0 & \begin{cases} \kk[-1] ~\text{if}~ a < c \leq b    \\ 0~\text{else} \end{cases} & \begin{cases} \kk ~\text{if}~ a < d \leq b    \\ 0~\text{else} \end{cases} & 0 & 0
\\

(-\infty, b) & \begin{cases} \kk[-1] ~\text{if}~  d < b    \\ 0~\text{else} \end{cases} & \begin{cases} \kk ~\text{if}~  c < b    \\ 0~\text{else} \end{cases} & \begin{cases} \kk ~\text{if}~  c < b \leq d    \\ 0~\text{else} \end{cases} & 0 & \begin{cases} \kk ~\text{if}~   b \leq d    \\ 0~\text{else} \end{cases}  & 0 & \kk & \begin{cases} \kk ~\text{if}~   c < b    \\ 0~\text{else} \end{cases} &  \kk \\

(a , +\infty) & \begin{cases} \kk[-1] ~\text{if}~  a < c    \\ 0~\text{else} \end{cases} & \begin{cases} \kk ~\text{if}~  a < d    \\ 0~\text{else} \end{cases} & 0 & \begin{cases} \kk ~\text{if}~  c \leq a < d    \\ 0~\text{else} \end{cases} & 0 & \begin{cases} \kk ~\text{if}~   a \leq c    \\ 0~\text{else} \end{cases} & \begin{cases} \kk ~\text{if}~  a < d    \\ 0~\text{else} \end{cases} & \kk & \kk \\

(-\infty, b] &  \begin{cases} \kk[-1] ~\text{if}~  c \leq b    \\ 0~\text{else} \end{cases} & \begin{cases} \kk ~\text{if}~  d \leq b    \\ 0~\text{else} \end{cases}  & 0 & \begin{cases} \kk[-1] ~\text{if}~  c \leq b < d   \\ 0~\text{else} \end{cases} & 0 & \begin{cases} \kk[-1] ~\text{if}~  c \leq b    \\ 0~\text{else} \end{cases} & \begin{cases} \kk ~\text{if}~  d \leq b    \\ 0~\text{else} \end{cases}  & 0 & 0 \\ 

[a, +\infty) & \begin{cases} \kk[-1] ~\text{if}~  a \leq d    \\ 0~\text{else} \end{cases} & \begin{cases} \kk ~\text{if}~  a \leq c    \\ 0~\text{else} \end{cases} & \begin{cases} \kk[-1] ~\text{if}~  c< a \leq d   \\ 0~\text{else} \end{cases} & 0  & \begin{cases} \kk[-1] ~\text{if}~   a \leq d   \\ 0~\text{else} \end{cases}  & 0 & 0 & \begin{cases} \kk ~\text{if}~   a \leq c   \\ 0~\text{else} \end{cases}  & 0 \\ 

\R & \kk[-1] & \kk  & 0 & 0 & 0 & 0 & \kk & \kk & \kk

\end{array}

\end{math}

\end{sideways}
\end{center}

 \end{prop}

\newpage
\normalsize

\begin{pf}
Since all computations work similarly, we shall only treat the case where $I = [a,b]$ and $J = (c,d) \not = \emptyset$. By proposition \ref{P:duality}: 

$$\Rr\Hom(\kk_I,\kk_J) \simeq \Rr \Gamma (\R ; \D (\D(\kk_J) \otimes \kk_I )). $$

We have $\D(\kk_J) \simeq \kk_{[c,d]}[1]$. Therefore, $\D (\D(\kk_J) \otimes \kk_I ) \simeq \D(\kk_{[c,d]\cap [a,b]}) [-1]$. If $[c,d]\cap [a,b] = \emptyset$, then $\Rr\Hom(\kk_I,\kk_J) = 0$. Otherwise, let $S =[c,d]\cap [a,b]$. Then $S$ is a closed interval. 

Let us first assume that $S$ has non empty interior $\text{Int}(S)$. Then $\D(\kk_{[c,d]\cap [a,b]}) [-1] \simeq \kk_{\text{Int}(S)}$, and from proposition \ref{P:sections}: $$ \Rr\Hom(\kk_I,\kk_J) \simeq \Rr \Gamma (\R ; \kk_{\text{Int}(S)}) \simeq \kk[-1].   $$

If $S$ has empty interior, since $S$ is a closed interval, there exists $t\in \R$ such that $S = \{t\}$. Therefore, from proposition \ref{P:duality}, we have $\D(\kk_{[c,d]\cap [a,b]}) [-1] \simeq \kk_{\{t\}} [-1]$. Consequently, we obtain from proposition \ref{P:sections}: 

$$ \Rr\Hom(\kk_I,\kk_J) \simeq \Rr \Gamma (\R ; \kk_{\{t\}}[-1]) \simeq \kk[-1].   $$

To sum up, we have proved that: $$ \Rr\Hom(\kk_I,\kk_J) \simeq \begin{cases} \kk[-1] ~\text{if}~ [a,b] \cap [c,d] \not = \emptyset  \\ 0~\text{else} \end{cases}. $$
\end{pf}

\subsection{The functor $- \star \K_\varepsilon$}

We start by recalling the following classical lemma about convolution of sheaves:

\begin{lem}[Exercice II.20 - \cite{Kash90}]\label{L:ConvolutionConvex}
Let $A,B\subset \V$ two closed subsets of the finite dimensional real vector space $\V$ (endowed with the topology inherited from any norm) satisfying: 
\begin{enumerate}
\item the map $s_{|A\times B} : A \times B \to \V$ is proper,
\item for any $x \in \V$, $s_{|A\times B}^{-1}(x)$ is contractible.
\end{enumerate}
 Then $\kk_A \star \kk_B \simeq \kk_{A+B}$ with $A+B = \{a+b \mid a\in A, b\in B\}$.
\end{lem}

Given $\varepsilon\geq 0$ and $I \subset \R$ an interval, we compute the convolution (Definition~\ref{D:Convolution}) $\kk_I \star \K_\varepsilon$ (see Equation \ref{eq:ballkernel}, Section \ref{S:Preliminaries}). The case where $I$ is closed is a direct consequence of lemma \ref{L:ConvolutionConvex}. We then deduce the other cases using distinguished triangles.

\begin{prop}\label{P:ComputationDerivedConvolution}

Let $\varepsilon \geq 0$, and $a\leq b$ in $\R \cup \{\pm \infty\}$, then : 
\begin{enumerate}

\item $\kk_{\R}\star \K_\varepsilon \simeq \kk_{\R} $,

\item $\kk_{[a,b]}\star \K_\varepsilon \simeq \kk_{[a-\varepsilon, b + \varepsilon]} $,

\item $\kk_{(a,b)}\star \K_\varepsilon \simeq 
\begin{cases}
\kk_{(a+\varepsilon,b-\varepsilon)} ~~\text{if} ~~\varepsilon < \frac{| b-a |}{2} \\
\kk_{[b - \varepsilon , a + \varepsilon]}[-1] ~~\text{if} ~~\varepsilon \geq \frac{| b-a |}{2}
\end{cases}$

\item $\kk_{(a,b]}\star \K_\varepsilon \simeq \kk_{(a+\varepsilon, b + \varepsilon]} $,

\item $\kk_{[a,b)}\star \K_\varepsilon \simeq \kk_{[a-\varepsilon, b - \varepsilon)} $.
\end{enumerate}

\end{prop}

\begin{pf}

We can obtain the computation for $\kk_{(a,b)}$ by using the distinguished triangle $\kk_{(a,b)} \longrightarrow  \kk_\R \longrightarrow \kk_{\R \backslash (a,b)} \stackrel{+1}{\longrightarrow}$, as  $\kk_{\R \backslash (a,b)}$ is the direct sum of one or two sheaves constant over closed intervals.
Similarly for the case of $\kk_{[a,b)}$, we can use the distinguished triangles $\kk_{[a,b)} \longrightarrow  \kk_{ [a,b]} \longrightarrow \kk_{\{b\}} \stackrel{+1}{\longrightarrow}$.
\end{pf}

\color{black}

\section{Structure of $\varepsilon$-interleavings}

In this section we investigate the structure of $\varepsilon$-interleavings between constructible sheaves over $\R$. We start by giving explicit conditions on the support and degree of two indecomposables sheaves $\kk_I[-i]$ and $\kk_J[-j]$ to be $\varepsilon$-interleaved. Then, we introduce for any constructible sheaf $F$ on $\R$ its CLR decomposition (Proposition \ref{P:CLRdecompositionexists}), which expresses $F$ as a direct sum of three sheaves $F_C$, $F_L$ and $F_R$ whose interval decomposition have specific properties. We further show that the CLR decomposition decomposes the notion of interleaving in the following sense: $F$ and $G$ are $\varepsilon$-interleaved if and only if $F_C$ and $G_C$, $F_L$ and $G_L$, $F_R$ and $G_R$ are (Theorem \ref{T:ConvolutiondistibutesoverParts}).
\color{black}

\subsection{Characterization of $\varepsilon$-interleavings between indecomposable sheaves}

For  any interval $I$ and real number $\varepsilon \geq 0$, we will write $I^\varepsilon = \mathop{\cup}\limits_{x \in I} \text{B}(x,\varepsilon)$ where $\text{B}(x,\varepsilon)$ is the euclidean closed ball centered at $x$ with radius $\varepsilon$. Moreover if $I = (a,b)$ with $a,b \in \R$, and $\varepsilon < \frac{b-a}{2} = \frac{\text{diam}(I)}{2}$, define $I^{-\varepsilon} = (a+\varepsilon , b-\varepsilon)$. If $I$ is bounded, we write $\text{cent}(I)$ for its center, that is $(b+a)/2$ where $a, b$ are the boundary points of $I$.

The following proposition describes the condition for sheaves constant on open/closed intervals to be $\varepsilon$-interleaved.  
\begin{prop}[closed/open]\label{P:closedOpen}
Let $S,T$ (resp. $U,V$) be two non-empty closed intervals (resp. non-empty open intervals), and $\varepsilon > 0$. Then:

\begin{enumerate}
\item $\kk_S \sim_\varepsilon \kk_T $ if and only if $ S \subset T^\varepsilon$ and $T \subset S^\varepsilon$,

\item $\kk_U \sim_\varepsilon \kk_V $ if and only if $ U \subset V^\varepsilon$ and $V \subset U^\varepsilon$.

\item Assuming that $S$ and $U$ are bounded,  $\kk_S \sim_\varepsilon \kk_U[-1] $ if and only if $ \varepsilon \geq \frac{\text{diam}(U)}{2}$ and $S \subset  \left[\text{cent}(U) - (\varepsilon - \frac{\text{diam}(U)}{2}), \text{cent}(U) + (\varepsilon - \frac{\text{diam}(U)}{2})\right]$.
\end{enumerate}
\end{prop}

\begin{pf}

\begin{enumerate}
\item Consider $f : \kk_S\star \K_\varepsilon \to \kk_T$ and $g : \kk_T\star \K_\varepsilon \to \kk_S$ the data of an $\varepsilon$-interleaving. Then $f$ and $g$ are in particular not zero since 
$$\G(\R, \phi_{\kk_S,2\varepsilon}) = \G(\R, g) \circ \G(\R, f\star \K_\varepsilon) $$
is an isomorphism between $\G(\R, \kk_S \star K_{2\varepsilon})$ and $\G(\R, \kk_T)$ which are non zero. 
 Remark that $\kk_S\star \K_\varepsilon \simeq \kk_{S^\varepsilon}$ and $\kk_T\star \K_\varepsilon \simeq \kk_{T^\varepsilon}$ by Proposition~\ref{P:ComputationDerivedConvolution}.  From our computations of morphisms (Proposition~\ref{P:TableHomNonDerived}), we have necessarily  $S \subset T^\varepsilon$ and $T \subset S^\varepsilon$. Conversely, if $S \subset T^\varepsilon$ and $T \subset S^\varepsilon$, it is easy to build an $\varepsilon$-interleaving.

\item Consider $f : \kk_U\star \K_\varepsilon \to \kk_V$ and $g : \kk_V\star \K_\varepsilon \to \kk_U$ the data of an $\varepsilon$-interleaving. For the same reason as above, $f$ and $g$ are not zero. Hence, $f\star K_{-\varepsilon} : \kk_U \to \kk_V \star K_{-\varepsilon}  $ is not zero. As $\kk_V \star K_{-\varepsilon} \simeq \kk_{V^\varepsilon} $, by our computations of morphisms (Proposition~\ref{P:TableHomNonDerived}) we get that $U \subset V^\varepsilon$. Similarly we have  $V \subset U^\varepsilon$.

Conversely if we assume $U \subset V^\varepsilon$ and $U \subset V^\varepsilon$, it is easy to construct an $\varepsilon$-interleaving.

\item Let $f : \kk_S\star \K_\varepsilon \to \kk_U[-1]$, $g : \kk_U[-1]\star \K_\varepsilon \to \kk_S$ be the data of an $\varepsilon$-interleaving. For the same reason as above, $f$ and $g$ are not zero. Suppose $\varepsilon < \frac{\text{diam}(U)}{2} $, then Proposition~\ref{P:ComputationDerivedConvolution} implies that $\kk_U[-1]\star \K_\varepsilon \simeq \kk_{U^{-\varepsilon}}[-1] $, hence the fact that $g$ is not zero is absurd.

Therefore we  have  $ \varepsilon \geq \frac{\text{diam}(U)}{2}$, and $$\kk_U[-1]\star \K_\varepsilon \simeq \kk_{[\text{cent}(U) - (\varepsilon - \frac{\text{diam}(U)}{2}), \text{cent}(U) + (\varepsilon - \frac{\text{diam}(U)}{2})]}.$$ Hence the existence of $g$ implies that $S \subset  \left[\text{cent}(U) - (\varepsilon - \frac{\text{diam}(U)}{2}), \text{cent}(U) + (\varepsilon - \frac{\text{diam}(U)}{2})\right]$. Also the existence of $f$ implies that $S^\varepsilon \cap U \not = \emptyset$, but this condition is weaker than the previous one. 

Conversely, if $ \varepsilon \geq \frac{\text{diam}(U)}{2}$ and $S \subset  \left[\text{cent}(U) - (\varepsilon - \frac{\text{diam}(U)}{2}), \text{cent}(U) + (\varepsilon - \frac{\text{diam}(U)}{2})\right]$, we can construct the desired morphisms (using Proposition~\ref{P:TableHomDerived}) and have to check that their composition (after applying $\varepsilon$ convolution to one of the two) is not zero, which can be obtained by taking stalks at any $x\in J$.

\end{enumerate}

\end{pf}

\begin{prop}[half-open]\label{P:Half-Open}
Let $I = [a,b)$ and $J = [c,d)$ with $a,c\in \R$ and $b,d\in \R\cup\{+\infty\}$, and $\varepsilon \geq 0$. Then $\kk_I \sim_\varepsilon \kk_J \iff \mid a-c \mid \leq \varepsilon $ and $ \mid b-d \mid \leq \varepsilon $.

Similarly for $I = (a,b]$ and $J = (c,d]$, $\kk_I \sim_\varepsilon \kk_J \iff \mid a-c \mid \leq \varepsilon $ and $ \mid b-d \mid \leq \varepsilon $.

\end{prop}

\begin{pf}
The proof works exactly the same as the open/closed case, that is Proposition~\ref{P:closedOpen}.
\end{pf}

\subsection{CLR Decomposition} 

In order to define a matching between graded barcodes, we have to distinguish between the topological nature of their support interval as the existence of shifted morphisms 
between them precisely depends on this nature (proposition \ref{P:TableHomDerived}).

\begin{defi}\label{D:CLRdecomposition}
Let $I\subset \R$ be an interval. 
\begin{enumerate}
\item $I$ is said to be an interval of type C if there exists $(a,b)\in \R^2$ such that $I = [a,b]$ or $I = (a,b)$.

\item $I$ is said to be an interval of type L if there exists $(a,b)\in \R^2$ such that $I = (a,b]$, $I = (-\infty, b]$ or $I = (a, + \infty)$.

\item  $I$ is said to be an interval of type R if there exists $(a,b)\in \R^2$ such that either $I = [a,b)$, $I = (-\infty, b)$, $I = [a,+\infty)$ or $I = \R$.

\end{enumerate}
\end{defi}

\begin{prop}\label{P:CLRdecompositionexists}
For $F \in \D^b_{\R c}(\kk_\R)$, there exists a decomposition, unique up to isomorphism,  $F \simeq F_C\oplus F_R \oplus F_L$ such that:
\begin{enumerate}
\item the cohomology objects of $F_C$ are direct sums of constant sheaves over intervals of type $C$,
\item the cohomology objects of $F_L$ are direct sums of constant sheaves over intervals of type $L$,
\item the cohomology objects of $F_R$ are direct sums of constant sheaves over intervals of type $R$.

\end{enumerate}

We will call $F_C$ (resp.  $F_L$, $F_R$) the {central} (resp. {left}, {right} ) part of $F$, and name this splitting the {CLR decomposition} of $F$.

\end{prop}

\begin{pf}
Observe that the types C,L,R do form a partition of the set of intervals of $\R$, and apply the decomposition and structure theorems from section~\ref{P:SheafKS}.
\end{pf}

\begin{defi}\label{D:CLRSheaf}
Let $F\in \D^b_{\R c}(\kk_\R)$. $F$ is said to be a {central sheaf} if $F\simeq F_C$. Similarly, $F$ is a {left} (resp. {right}) sheaf if $F \simeq F_L$ (resp. $F \simeq F_R$).
\end{defi}
We have the following easy properties.
\begin{prop}\label{P:CLRconvolution}
Let $F \in \D^b_{\R c}(\kk_\R)$, $\varepsilon \in \R $ and $\alpha \in \{C,L,R\}$. Then $F$ is of type $\alpha$ if and only if $F \star \K_\varepsilon$ is of type $\alpha$.
\end{prop}

\begin{pf}
It is sufficient to prove the statement for sheaves of the form $\kk_I[i]$ for $I \subset \R$ an interval and $i\in \Z$, which is a direct consequence of proposition \ref{P:ComputationDerivedConvolution}.
\end{pf}

\begin{prop}\label{P:CLRmorphism}
Let $F$ and $G$ be two sheaves of type $\alpha \in \{C,L,R\}$. Then any morphism $F \to G$ that factorizes through a sheaf $H$ as $F \to H \to G$, with $H$ of type $\beta \in \{C,L,R\}\backslash\{\alpha\}$, is necessarily zero. 
\end{prop}

\begin{pf}
It is sufficient to prove the statement for $F= \kk_I[i]$, $G = \kk_J[j]$ and $H = \kk_L[l]$, with $I$ and $J$ two intervals of type $\alpha$, $L$ an interval of type $\beta$, and $i,j,l \in \Z$. This is then a direct consequence of our computations of morphisms in $\D^b_{\R c }(\kk_\R)$ (proposition \ref{P:TableHomDerived}). 
\end{pf}

The CLR decomposition is compatible with the relation of being $\varepsilon$-interleaved in the following sense.

\begin{thm} \label{T:ConvolutiondistibutesoverParts}
Let $F,G \in \D^b_{\R c}(\kk_\R))$ and $\varepsilon \geq 0$, then the following holds : 

$$F \sim_\varepsilon G \iff \begin{cases}
F_C \sim_\varepsilon G_C  \\F_L \sim_\varepsilon G_L \\   F_R \sim_\varepsilon G_R \end{cases}. $$

\end{thm}

\begin{pf}
The right to left implication is an immediate consequence of the additivity of the convolution functor. 

We choose two isomorphisms, that will remain the same through all the proof: $$ F \simeq F_C \oplus F_L \oplus F_R ~~~\text{and}~~~G \simeq G_C \oplus G_L \oplus G_R.$$

Then by proposition \ref{P:CLRconvolution}, with $ \eta \in \R$, applying the functor $- \star \K_\eta$ to the above isomorphisms gives us the CLR decompositions of $F \star \K_\eta$ and $G \star \K_\eta$ in terms of those of $F$ and $G$. Now let us consider the data of an $\varepsilon$-interleaving between $F$ and $G$, that is, two morphisms $F \star \K_\varepsilon \stackrel{f}{\longrightarrow} G$ and $G \star \K_\varepsilon \stackrel{g}{\longrightarrow} F $ such that $f \star \K_\varepsilon \circ g : G \star K_{2\varepsilon} \longrightarrow G$ is the smoothing morphism $\phi_{G, 2 \varepsilon}$ (see definition \ref{D:smoothingmorphism}) and similarly  $g \star \K_\varepsilon \circ f : F \star K_{2\varepsilon} \longrightarrow F$ equals $\phi_{F, 2 \varepsilon}$. For $\alpha , \beta \in \{C,L,R\}$, we denote by $f_{\alpha,\beta}$ the composition 

$$f_{\alpha,\beta} = F_\alpha \star \K_\varepsilon \longrightarrow (F_C \oplus F_L \oplus F_R )\star \K_\varepsilon \simeq  F \star \K_\varepsilon \stackrel{f}{\longrightarrow} G \simeq G_C \oplus G_L \oplus G_R \longrightarrow G_\beta. $$

We denote $f_{\alpha,\alpha}$ by $f_\alpha$, and use similar notations for $g$ and the smoothing morphisms of $F$ and $G$. By proposition \ref{P:CLRmorphism}, we have: 

$$g \star \K_\varepsilon \circ f = g_C \star \K_\varepsilon \circ f_C + g_L \star \K_\varepsilon \circ f_L + g_R \star \K_\varepsilon \circ f_R,  $$
$$f \star \K_\varepsilon \circ g = f_C \star \K_\varepsilon \circ g_C + f_L \star \K_\varepsilon \circ g_L + f_R \star \K_\varepsilon \circ g_R. $$

Therefore, for $\alpha \in \{C,L,R\}$, projecting the above equations onto the summands of type $\alpha$ gives:
$$g_\alpha \star \K_\varepsilon \circ f_\alpha = \phi_{F_\alpha, 2 \varepsilon} ~~~\text{and}~~~ f_\alpha \star \K_\varepsilon \circ  g_\alpha = \phi_{G_\alpha, 2 \varepsilon},$$

since $(\phi_{F, 2 \varepsilon})_{\alpha} = \phi_{F_\alpha, 2 \varepsilon}$ and $(\phi_{G, 2 \varepsilon})_{\alpha} = \phi_{G_\alpha, 2 \varepsilon}$. Consequently, we deduce that $F_\alpha \sim_\varepsilon G_\alpha$.
\end{pf}

\section{Isometry theorem and graded barcodes }

This section presents the proof of the isometry theorem problem associated to the Krull-Schmidt metric category $(\D^b_{\R c}(\kk_\R), d_C)$ (see section \ref{S:IsomProblem}). The inequality $d_C \leq d_B$ is an easy consequence of the additivity of the convolution functor. 

In order to prove the reverse inequality, we prove that an $\varepsilon$-interleaving between 
two sheaves induces a $\varepsilon$-matching between their graded-barcodes. To do so, we construct the matching according to the CLR decomposition. We reduce the construction of the matching between the left and right parts to the well-known case of persistence modules with one parameter. To this end, we first prove that interleavings between right (resp. left) parts of two sheaves happen degree-wise at the level of their cohomology objects. This enables us to define functors $\Psi^j_R$, that send the $j$-th cohomology of the right part of a sheaf to a one parameter persistence module. We prove that $\Psi^j_R$ are barcode preserving, and send interleavings of sheaves to interleavings of persistence modules.

\subsection{The easy inequality}

We start by proving the easy direction of the inequality.

\begin{lem}\label{L:ConvsmallerBottleneck}
Let $F$ and $G$ two objects of $\D^b_{\R c}(\kk_\R)$, then: $$ d_C(F,G) \leq d_B(\B(F),\B(G)).$$

\end{lem}

\begin{pf}
If $d_B(\B(F),\B(G)) = +\infty$, then the inequality holds. Let us now assume that  $d_B(\B(F),\B(G)) < +\infty$.
Let $\delta \in \{\varepsilon \geq 0 \mid \text{there exists a } \varepsilon  \text{-matching between } \B(F) \text{ and } \B(G) \}$. Then for any $\eta > \delta$ there exists two subsets $\mathcal{X} \subset \B(F)$ and $\mathcal{Y} \subset \B(G)$ and a bijection $\sigma : \mathcal{X} \to \mathcal{Y}$ such that for all $I^i\in \mathcal{X}$, with $\sigma(I^i) = J^j$ there exists an $\eta$-interleaving between $\kk_I[-i]$ and $\kk_{J}[-j]$ given by the two morphisms :

$$f_{I^i} : \kk_I[-i]\star \K_\eta \to \kk_{J}[-j] ~~~\text{and}~~~g_{I^i} : \kk_J[-j]\star \K_\eta \to \kk_{I}[-i], $$

and for all $K^k \in \B(F) \backslash \mathcal{X} \sqcup \B(G) \backslash  \mathcal{Y}  $, $\kk_K[-k] \sim_\eta 0 $.

Therefore, the morphisms:

$$f : F \star \K_\eta \longrightarrow \bigoplus_{I^i \in \mathcal{X}} \kk_I[-i] \star \K_\eta \stackrel{\oplus f_{I^i}}{\longrightarrow} \bigoplus_{J^j\in \mathcal{Y}} \kk_J[-j] \longrightarrow G $$

$$g : G \star \K_\eta \longrightarrow \bigoplus_{J^j \in \mathcal{Y}} \kk_J[-j] \star \K_\eta \stackrel{\oplus g_{I^i}}{\longrightarrow} \bigoplus_{I^i\in \mathcal{X}} \kk_I[-i] \longrightarrow F $$

form an $\eta$-interleaving between $F$ and $G$. Consequently, 

$$ \{\varepsilon \geq 0 \mid \text{there exists a } \varepsilon  \text{-matching between } \B(F) \text{ and } \B(G) \}$$ $$ \subset  \{\varepsilon \geq 0 \mid \text{there exists a } \varepsilon  \text{-isomorphism between } F \text{ and } G \}$$

which proves the lemma by taking the infimum of both sets.
\end{pf}

\subsection{The cases $F_R \leftrightarrow G_R$ and $F_L \leftrightarrow G_L$}

 In this section, we give a description of the $\varepsilon$-interleavings between the right parts of two complexes of sheaves. The proofs and statements for the left parts  are analogous.
 
\subsubsection{Construction of $\Psi^j_R$}
\begin{prop}
Let $F,G \in \D^b_{\R c}(\kk_\R)$ and $\varepsilon \geq 0$ with right parts $F_R$ and $G_R$. The following holds : 

$$F_R \sim_\varepsilon G_R \iff \forall j \in \Z,~ \Ho^j(F_R) \sim_\varepsilon \Ho^j(G_R). $$

\end{prop}

\begin{pf}
The right to left implication is clear, so let us consider an $\varepsilon$-interleaving given by $F_R\star \K_\varepsilon \stackrel{f}{\longrightarrow} G_R$ and $G_R\star \K_\varepsilon \stackrel{g}{\longrightarrow} F_R$. Let $j \in \Z$ and pick $\kk_I$ a direct summand of $\Ho^j(F_R)$ ($I$ is a half-open interval of the type $[a,b)$). We consider again the composition : 

$$ \kk_I[-j]\star K_{2\varepsilon} \stackrel{i^{F_R}_I}{\longrightarrow} F_R \star K_{2\varepsilon}  \stackrel{f \star \K_\varepsilon}{\longrightarrow} G_R\star \K_\varepsilon \stackrel{g}{\longrightarrow} F_R  \stackrel{p^{F_R}_I}{\longrightarrow} \kk_I[-j]$$

From our computations of derived morphisms (Proposition~\ref{P:TableHomDerived}), this is equal to :  

$$ \kk_I[-j]\star K_{2\varepsilon} \stackrel{i^{F_R}_I}{\longrightarrow} F_R \star K_{2\varepsilon}  \stackrel{f \star \K_\varepsilon}{\longrightarrow}  \Ho^j(G_R\star \K_\varepsilon)[-j] \oplus \Ho^{j+1}(G_R\star \K_\varepsilon)[-j-1]  \stackrel{g}{\longrightarrow} F_R  \stackrel{p^{F_R}_I}{\longrightarrow} \kk_I[-j].$$

We obtain using our computations of convolution (Proposition~\ref{P:ComputationDerivedConvolution})  that, since $G_R$ has only half-open intervals in the decomposition of its cohomology objects, $\Ho^j(G_R\star \K_\varepsilon)[-j] \simeq \Ho^j(G_R)[-j] \star \K_\varepsilon$ and  $\Ho^{j+1}(G_R\star \K_\varepsilon)[-j-1] \simeq \Ho^{j+1}(G_R)[-j-1] \star \K_\varepsilon$. 

It follows again from  Proposition~\ref{P:TableHomDerived} that any morphism of $\kk_I[-j]\star \K_\varepsilon \to \kk_I[-j] $ that factors through a complex concentrated in degree $j+1$ must be zero. 

Finally, the first composition is thus equal to 

$$ \kk_I[-j]\star K_{2\varepsilon} \stackrel{i^{F_R}_I}{\longrightarrow} F_R \star K_{2\varepsilon}  \stackrel{f \star \K_\varepsilon}{\longrightarrow}  \Ho^j(G_R)[-j] \star \K_\varepsilon \stackrel{g}{\longrightarrow} F_R  \stackrel{p^{F_R}_I}{\longrightarrow} \kk_I[-j].$$

As this is true for any summand of $\Ho^j(F_R)$ we get that the composition :

$$ \Ho^j(F_R)[-j] \star K_{2\varepsilon} \stackrel{}{\longrightarrow} F_R \star K_{2\varepsilon}  \stackrel{f \star \K_\varepsilon}{\longrightarrow} G_R\star \K_\varepsilon \stackrel{g}{\longrightarrow} F_R  \stackrel{}{\longrightarrow} \Ho^j(F_R)[-j]$$
is equal to the composition 
$$ \Ho^j(F_R)[-j]\star K_{2\varepsilon} \stackrel{}{\longrightarrow} F_R \star K_{2\varepsilon}  \stackrel{}{\longrightarrow} \Ho^j(G_R)[-j]\star \K_\varepsilon \stackrel{}{\longrightarrow} F_R  \stackrel{}{\longrightarrow} \Ho^j(F_R)[-j].$$
This gives the first part of the $\varepsilon$-interleaving. We get the second one by intertwining the roles of $F_R$ and $G_R$.
\end{pf}

The result above shows that when one wants to understand a morphism between the right parts of two sheaves, it is sufficient to  understand it at the level of each of their cohomology objects, degree wise. We will show that the behavior of $\varepsilon$-interleavings between sheaves with cohomologies concentrated in degree $j \in \Z$ decomposing into direct summands of type $R$, is essentially the same as looking at $\varepsilon$-interleavings in the opposite category of one-parameter persistence modules, which is well understood. We quickly introduce all the necessary definitions and results needed, but we refer to \cite{ChazSilv16} for a detailed exposition about the isometry theorem for one-parameter persistence.

We denote $\Pe_f(\kk^\R)$ the category of pointwise finite dimensional persistence modules over $\R$, that is, the category of functors $M : (\R, \leq) \longrightarrow \text{Mod}_f(\kk)$ where $\text{Mod}_f(\kk)$ is the category of finite dimensional vector spaces over the field $\kk$. 

There is a notion of $\varepsilon$-interleaving (for $\varepsilon \geq 0$ ) in this context based on the shift functor $\cdot[\varepsilon]$ defined as $M[\varepsilon](s) = M(s +\varepsilon)$ and $M[\varepsilon](s\leq t)= M(s+\varepsilon \leq t+\varepsilon)$ for $s\leq t$ two real number. There is also a canonical natural transformation $s^M_\varepsilon : M \longrightarrow M[\varepsilon] $. We will say that $M$ and $N$ in $\Pe_f(\kk^\R)$ are $\varepsilon$-interleaved if there exists two morphisms $f : M\longrightarrow N[\varepsilon]$ and $g : N \longrightarrow M[\varepsilon]$ such that $g[\varepsilon] \circ f = s^M_{2 \varepsilon}$ and $f[\varepsilon] \circ g = s^N_{2 \varepsilon}$.

The pseudo-distance induced on  $\Pe_f(\kk^\R)$ by:
$$d_I(M,N) := \inf \{\varepsilon\geq 0 \mid M ~\text{and}~N ~\text{are}~\varepsilon-\text{interleaved}\}$$
is called the interleaving distance, and was first introduced in \cite{CSGG09}.
In this text, for $I \subset \R$ an interval, we will write $\kk^I$ the object of $\Pe_f(\kk^\R)$ defined, for $s\leq t$, by: 
$$ \kk^I(s) = \begin{cases} \kk ~\text{if}~s\in I
\\ 0~\text{else}
\end{cases} ~~~ \kk^I(s\leq t) = \begin{cases} \text{id}_\kk ~\text{if}~s,t\in I
\\ 0~\text{else}
\end{cases}. $$

Let $\D^j_R$ be the full sub-category of $\D^b_{\R c}(\kk_\R)$ of complexes of sheaves $F$ such that $F \simeq F_R$ and $\Ho^i(F) = 0$ for $i \not = j$.

\begin{prop}\label{P:ExistenceofPsi}
There exists a functor $\Psi^j_R : \D^j_R \longrightarrow \Pe(\kk_\R)^{op}$ such that : 

\begin{enumerate}
\item for $F \in \D^j_R$ such that $\Ho^j(F)\simeq \oplus_{I \in B}\kk_I$, we have $\Psi^j_R(F) = \oplus_{I \in B}\kk^I$, 

\item $\Psi^j_R$ is fully faithful,

\item for $\varepsilon \geq 0$ and $F \in \D^j_R$, $\Psi^j_R(F \star \K_\varepsilon ) = \Psi^j_R(F)[\varepsilon]$ and  $\Psi^j_R(\phi_{F,\varepsilon}) = s_\varepsilon^{\Psi^j_R(F)}$,

\item $\Psi^j_R$ is isometric with respect to $d_C(\cdot,\cdot)$ and $d_I(\cdot,\cdot)$.
\end{enumerate}

\end{prop}

\begin{pf}
This is a combination of the computations of morphisms and convolution (Propositions~\ref{P:TableHomNonDerived}, \ref{P:ComputationDerivedConvolution}, \ref{P:CLRdecompositionexists}), together with the observation that for $I,J$ two intervals of type $R$ and $j \in \Z$, then we have the functorial isomorphisms : $$\Hom_{\D^b_{\R c}(\kk_\R)}(\kk_I[-j],\kk_J[-j]) \simeq \Hom_{\Mod(\kk_\R)}(\kk_I,\kk_J) \simeq \Hom_{\Pe(\kk_\R)}(\kk^J,\kk^I). $$ 

\end{pf}

\subsubsection{Matching of the right parts}

\begin{thm}[Matching of right parts]
Let $F,G \in \D^b_{\R c}(\kk_\R)$ be $\varepsilon$-interleaved with respect to the morphisms
$F \star \K_\varepsilon \stackrel{f}{\longrightarrow} G$ and $G \star \K_\varepsilon \stackrel{g}{\longrightarrow} F $. 
Let $j\in \Z$. Then there exists an $\varepsilon$-matching $\sigma^j_R : \B^j(F_R) \not \longrightarrow \B^j(G_R) $ (see definition \ref{D:matching}).
\end{thm}

\begin{pf}
Observe that $\Psi^j_R(F_R)$ (resp. $\Psi^j_R(G_R)$) is a persistence module with the same barcode than $\Ho^j(F_R)$ (resp. $\Ho^j(G_R)$). Also, from proposition 5.3, $\Psi^j_R(F_R)$ and $\Psi^j_R(G_R)$ are $\varepsilon$-interleaved as persistence modules. Hence, we can apply the isometry theorem for pointwise finite dimensional persistence modules \cite[Theorem 4.11]{ChazSilv16} to $\Psi^j_R(F_R)$ and $\Psi^j_R(G_R)$ and deduce the existence of a $\varepsilon$-matching of barcodes of persistence modules between   $\B^j(F_R)$ and $\B^j(G_R)$.
This matching is  what we ask for $\sigma^j_R$ by Proposition~\ref{P:Half-Open}.
\end{pf}

\subsection{The case $F_C\leftrightarrow G_C$}

In this section, we construct the $\varepsilon$-matching between the central parts of two sheaves, assuming they are $\varepsilon$-interleaved. Using ideas of Bjerkevik \cite[Section 4]{Bjer16}, we introduce a pre-order $\leq_\alpha$ on the set of graded-intervals of type C whose purpose is to prove the existence of the $\varepsilon$-matching using Hall's marriage  theorem. To do so, we must prove that given a finite list of interval in the barcode of one of the two sheaves, there exists, at least, the same number of intervals in the barcode of the second sheaf which are at distance less than $\varepsilon$ from an interval in the first list. 

We will show that ordering the graded-barcodes of the central sheaves according to $\leq_\alpha$ will actually lead to a very nice expression of the interleaving morphisms, allowing us, by a rank argument, to deduce that this condition is satisfied.

\subsubsection{Ordering graded-intervals of type C}
Recall that we defined a graded interval to be an interval $I$ together with an integer $j\in\Z$. It will be written $I^j$ henceforth. For $I$ of type C such that either $I = [a,b]$ or $I = (a,b)$ with $a,b\in \R$, define $\text{diam}(I) = b -a$ to be its diameter.

\begin{defi}\label{D:preorder}
The relation $\leq_\alpha $ on the set of graded intervals of type C is defined by : 

\begin{enumerate}
\item  For $R^i,T^j$ two closed intervals in degree $i$ and $j$ : $R^i \leq_\alpha T^j \iff i=j$ and $\text{diam}(T) \leq \text{diam}(R)$,
\item
for $U^i,V^j$ two open intervals in degree $i$ and $j$ : $U^i \leq_\alpha V^j \iff i=j$ and $\text{diam}(U) \leq \text{diam}(V)$,

\item for $R^i$ a closed interval in degree $i$, and $V^j$ an open interval in degree $j$ : $R^i \leq_\alpha U^j \iff i = j+1$.  
\end{enumerate}
\end{defi}

\begin{prop}\label{P:Preorder} The relation
$\leq_\alpha$ is a partial pre-order over the set of graded intervals of type C, that is, it is reflexive and transitive. Moreover, it is total if restricted to sets of graded intervals containing only, for a given $i\in \Z$, open intervals in degree $i$ and closed intervals in degree $i+1$. 
\end{prop}

The following is the analogous result in our setting to \cite[Lemma 4.6]{Bjer16}.

\begin{prop}
Let $I^i,J^i,S^l$ be three graded intervals of type $C$ and $\varepsilon\geq 0$ such that $I^i\leq_\alpha J^j$ and there exists two non-zero morphisms
$ \chi : \kk_S[-l] \star \K_\varepsilon \longrightarrow \kk_I[-i] ~~\text{and} ~~\xi : \kk_J[-j]\star \K_\varepsilon \longrightarrow \kk_S[-l]$.
Then either $\kk_S[-l] \sim_\varepsilon \kk_I[-i]$ or $\kk_S[-l] \sim_\varepsilon \kk_J[-j]$.
\end{prop}

\begin{pf}
By definition of the pre-order $\leq_\alpha$, we only have to investigate the three cases of the above definition~\ref{D:preorder} :

\begin{enumerate}
\item Let $i\in \Z$ and $R,T$ be two open intervals such that $R^i \leq_\alpha T^i$, that is, $\text{diam}(T) \leq \text{diam}(R)$. Let $S^l$ be a graded interval such that there exists some non-zero $\chi$ and $\xi$. Then $S$ must be a closed interval, and $l = i$. As a consequence, $R \subset S^\varepsilon$ and $S \subset T^\varepsilon$.

Assume that $\kk_R[-i] \not \sim_\varepsilon \kk_S[-i]$. Then, as $R \subset S^\varepsilon$, $S \not \subset R^\varepsilon$. So either $\min(S)< \min(R) - \varepsilon$, or $\max(S) > \max(R) + \varepsilon$. Assume the latter.

As $S \subset R^\varepsilon$, $\min(S)-\varepsilon < \min(R)$, we get subtracting the first inequality to this one : $\text{diam}(S)+\varepsilon > \text{diam}(R) + \varepsilon$. Hence $S <_\alpha R$. 
We get the same thing assuming $\min(S)< \min(R) - \varepsilon$.

Moreover, one can prove this way that $\kk_T[-i]\not \sim_\varepsilon \kk_S[-i]$ implies $S <_\alpha R$.

As we assumed $R^i \leq_\alpha T^i$, one has $\kk_S[-l] \sim_\varepsilon \kk_I[-i]$ or $\kk_S[-l] \sim_\varepsilon \kk_J[-j]$.

\item The proof for $U^i,V^i$ where $U$ and $V$ are open intervals is similar. 

\item Let $R^i$ a closed interval, $V^j$ an open interval, with $i = j+1$. Let $S^l$ be a graded interval and $\varepsilon$ such that there exists $\chi$ and $\xi$ such as in the proposition. Then $S$ must be an open interval and $l = j$. By the existence of $\chi$, we have that $\varepsilon \geq \frac{\text{diam}(U)}{2}$ and $R \subset [\text{cent}(U) - (\varepsilon - \frac{\text{diam}(U)}{2}), \text{cent}(U) + (\varepsilon - \frac{\text{diam}(U)}{2})]$, which, according to our characterization of $\varepsilon$-interleaving between indecomposable sheaves (proposition \ref{P:closedOpen}), is equivalent to $\kk_R[-j-1] \sim_\varepsilon \kk_S[-j]$. 

\end{enumerate}

\end{pf}
\subsubsection{Induced matching}

We now have the ingredients to prove the theorem. We start by introducing a sign notation. Given two intervals $I$ and $J$ of type C, we define: 

$$\delta(I,J) = \begin{cases}
0 ~\text{if $I$ and $J$ are both closed or both open,} \\
1~\text{if $I$ is open and $J$ is closed}, \\
-1 ~\text{if $I$ is closed and $J$ is open}.
\end{cases} $$

\begin{thm}[Matching of central parts]\label{T:matchingcentralparts}

Let $F_C$ and $G_C$ be two central sheaves (definition \ref{D:CLRSheaf}), and $\varepsilon \geq 0$ be such that $F_C$ and $G_C$ are $\varepsilon$-interleaved with respect to maps $F_C \star \K_\varepsilon \stackrel{f}{\longrightarrow} G_C $ and $G_C \star \K_\varepsilon \stackrel{g}{\longrightarrow} F_C $. Then,  there exists a bijection  
$$\sigma_C :  {\B}(F_C) \longrightarrow {\B}(G_C)$$
such that, for $I^j \in \mathbb{B}^j_C $, with $J = \sigma_C(I)$, we have $J^{j + \delta(I,J)} \in \mathbb{B}^{j + \delta(I,J)}(G_C)$ and $\kk_I \sim_\varepsilon \kk_{J}[-\delta(I,J)]$.

\end{thm}

Our proof will use a generalization of Hall's marriage theorem to the case of countable sets. For a reference, see for instance \cite{Pode76}. 

\begin{thm}[Hall]
Let $X$ and $Y$ be two countable sets, let $\mathcal{P}(Y)$ be the set of subsets of $Y$ and $M : X \to \mathcal{P}(Y)$. Then the following are equivalent : 

\begin{enumerate}
\item there exists an injective map $m : X \to Y$ satisfying $m(x) \in M(x)$ for every $x\in X$;

\item for every finite subset $A\subset X$, $|A|\leq |\cup_{x\in A}M(x)|$. Where $|A|$ is the cardinality of $A$.
\end{enumerate}

\end{thm}

We let $F_C$ and $G_C$ be two central sheaves. We set two isomorphisms : 
$$ F_C \simeq  \bigoplus_{I^j \in \mathbb{B}(F_C)} \kk_I[-j] ~~~\mbox{and} ~~~ G_C \simeq  \bigoplus_{I^j \in \mathbb{B}(G_C)} \kk_I[-j].  $$ 

For any morphism $f : F_C \to G_C$, given $I^i\in \B(F_C)$ and $J^j\in \B(G_C)$, we will write : 
$$f_{I^i,J^j} = \kk_I[-i] \longrightarrow F_C \stackrel{f}{\longrightarrow} G_C \longrightarrow \kk_J[-j]. $$

Similarly for $A \subset \B(F) $, let $f_{|A}$ be the composition : 

$$\bigoplus_{I^i \in A} \kk_I[-i] \longrightarrow F_C \stackrel{f}{\longrightarrow} G_C.  $$

We now assume that $F_C$ and $G_C$ are $\varepsilon$-interleaved with respect to $F_C \star \K_\varepsilon \stackrel{f}{\longrightarrow} G_C $ and $G_C \star \K_\varepsilon \stackrel{g}{\longrightarrow} F_C $. For $I^i\in \B(F_C) $ and $J^j\in \B^j(G_C)$, we deduce from our computations of propositions  \ref{P:TableHomDerived} and \ref{P:closedOpen} that: $$(f\star \K_\varepsilon)_{I^i,J^j} \circ g_{J^j,I^i} \not = 0 \text{ implies either : }
\begin{cases} I,J ~\text{are closed and}~ i=j, \\
I \text{ is open}, J \text{ is closed and } j = i+1, \\
J \text{ is open}, I \text{ is closed and } i = j+1.
\end{cases} $$

\begin{pf}[matching of central parts]
Our strategy is to adapt Bjerkevik's proof of \cite[Theorem 4.2]{Bjer16} to our setting. The pre-order $\leq_\alpha$ we have defined has exactly the same properties as the one defined in his proof. 

To define $\sigma_C$, we will apply Hall's theorem. From the local finiteness properties follows the fact that the graded-barcodes of $F_C$ and $G_C$ are countable. We here consider multi-sets as sets, to make the proof easier to understand. Nevertheless, it would not be difficult to write the proof properly using multi-sets. Let $M : {\B}_C(F_C) \to \mathcal{P}({\B}_C(G_C))$ defined by :  

$$M(I^i) = \{J^j \in {\B}(G_C) \mid \kk_I[-i] \sim_\varepsilon \kk_{J}[-j-\delta(I,J)]\} $$
for $I^i \in {\B}(F_C)$.

We define the following partitions : ${\B}(F_C) = \mathop{\sqcup}\limits_{i \in \Z} \Sigma^{i}_{F_C}$ and ${\B}(G_C) = \mathop{\sqcup}\limits_{i \in \Z} \Sigma^{i}_{G_C}$ where, 

$$\Sigma^{i}_{F_C} = \{J^j \in \B(F_C) \mid \text{$J$ is open and $j=i$ or $J$ is closed and $j = i+1$ }\}, $$
$$\Sigma^{i}_{G_C} = \{J^j \in \B(G_C) \mid \text{$J$ is open and $j=i$ or $J$ is closed and $j = i+1$ }\}. $$

We will define $\sigma_C$ according to these partitions, that is, we will construct some bijections $\sigma_C^i :\Sigma^{i}_{F_C} \longrightarrow \Sigma^{i}_{G_C} $ for all $i\in \Z$ and set $\sigma_C = \sqcup_{i \in \Z} \sigma_C^i$.

Let $i\in \Z$,  $A$ be a finite subset of $\Sigma^{i}_{F_C}$ and $M(A) = \cup_{I^i\in A} M(I^i)$. To apply Hall's theorem and deduce the existence of $\sigma_C$, we need to prove that $|A|\leq |M(A)|$.

By proposition \ref{P:Preorder}, $\leq_\alpha$ is a total pre-order on $A$. Hence, with $r = |A|$, there exists an enumeration $A = \{I_1^{i_1},...,I_r^{i_r}\}$, where $i_l = n$ if $I_l$ is an open interval and $i_l = n +1$ if $I_l$ is a closed interval, such that for $1\leq i\leq j \leq r$ we have $ I_i \leq_\alpha I_j$.

 We have by assumption $g \circ (f \star \K_\varepsilon) = \phi_{F,2\varepsilon}$ (see definition \ref{D:inteleavingconvolution}), also, the additivity of the convolution functor implies the following equality for $I_l^{i_l} \in A$ : 

$$ \phi_{\kk_I[-i_l],2\varepsilon} = \kk_I[-i_l]\star \K_\varepsilon \longrightarrow F\star \K_\varepsilon  \stackrel{\phi_{F,2\varepsilon}}{\longrightarrow} F \to \kk_I[-i_l]. $$ 

Therefore : 

\begin{align*}
\phi_{\kk_I[-i_l],2\varepsilon} &=  \sum_{J^j \in \B(G) } g_{J^j,I^{i_l}_l} \circ (f \star \K_\varepsilon)_{I^{i_l}_l,J^j}  \\
 &= \sum_{J^j\in {\B^l}(G)}  g_{J^j,I^{i_l}_l}\circ  \left [f_{I^{i_l}_l,J^j}  \star \K_\varepsilon \right].
\end{align*}

Now observe that if $g_{J^j,I^{i_l}_l}\circ  \left [f_{I^{i_l}_l,J^j}  \star \K_\varepsilon \right ] \not = 0 $ then $\kk_{I_l}[-i_l] \sim_\varepsilon \kk_J[-j]$, hence : 

\begin{align*}
\phi_{\kk_I[-i_l],2\varepsilon} &=  \sum_{J^j \in M(A) } g_{J^j,I^{i_l}_l} \circ (f \star \K_\varepsilon)_{I^{i_l}_l,J^j} 
\end{align*}

Similarly for $I_m \not = I_{m'} $ in $A$, $$0=\sum_{J^j\in {\B^l}(G)}  g_{J^j,I_m^{i_m}}\circ  \left [f_{I_{m'}^{i_m'},J^j}  \star \K_\varepsilon \right].$$

Hence if $m < m'$ and $g_{J^j,I_m^{i_m}}\circ  \left [f_{I_{m'}^{i_m'},J^j}  \star \K_\varepsilon \right ] \not = 0$, then $\kk_J[-j]$ is $\varepsilon$-interleaved with either $\kk_{I_m}[-m]$ or $\kk_{I_{m'}}[-m']$. Therefore : 
$$0=\sum_{J^j\in M(A)}  g_{J^j,I_m^{i_m}}\circ  \left [f_{I_{m'}^{i_m'},J^j}  \star \K_\varepsilon \right ]. $$

For $m> m'$, we can't say anything about the value of $\sum_{J^j\in M(A)}  g_{J^j,I_m^{i_m}}\circ \left [ f_{I_{m'}^{i_m'},J^j}  \star \K_\varepsilon \right ]$.

Writing those equalities in matrix form, we get :

\begin{bigcenter}
\begin{minipage}{20cm}
\small 
$$\left(\begin{array}{ccc}g_{J^1,I_{1}^{i_1}}  & \dots & g_{J^1,I_{r}^{i_r}}   \\\vdots & \ddots & \vdots \\g_{J^s,I_{1}^{i_1}}   & \dots & g_{J^s,I_{r}^{i_r}}  \end{array}\right)  \left(\begin{array}{ccc}f_{I_{1}^{i_1},J^1}  \star \K_\varepsilon & \dots & f_{I_{r}^{i_r},J^1}  \star \K_\varepsilon \\\vdots & \ddots & \vdots \\f_{I_{1}^{i_1},J^s}  \star \K_\varepsilon & \dots & f_{I_{r}^{i_r},J^s}  \star \K_\varepsilon\end{array}\right) =\left(\begin{array}{cccc}\phi_{I_1^{i_1},2\varepsilon} & ? & ? & ? \\0 & \phi_{I_2^{i_2},2\varepsilon} & ? & ? \\\vdots & \vdots & \ddots & ? \\0 & 0 & \dots & \phi_{I_r^{i_r},2\varepsilon}\end{array}\right) $$

\end{minipage}
\end{bigcenter}
\normalsize

Now recall that $\G(\R,-)$ is an additive functor. Hence, applying $\G(\R,-)$ to the above equality, we get : 

\small 

\begin{align*}
\left(\begin{array}{ccc}\G(\R,g_{J^1,I_{1}^{i_1}})  & \dots & \G(\R,g_{J^1,I_{r}^{i_r}})   \\ \vdots & \ddots & \vdots \\\G(\R,g_{J^s,I_{1}^{i_1}})   & \dots & \G(\R,g_{J^s,I_{r}^{i_r}})  \end{array}\right) & \left(\begin{array}{ccc}\G(\R,f_{I_{1}^{i_1},J^1}  \star \K_\varepsilon) & \dots & \G(\R,f_{I_{r}^{i_r},J^1}  \star \K_\varepsilon )\\\vdots & \ddots & \vdots \\\G(\R,f_{I_{1}^{i_1},J^s}  \star \K_\varepsilon) & \dots & \G(\R,f_{I_{r}^{i_r},J^s}  \star \K_\varepsilon)\end{array}\right)  \\ &=
 \left(\begin{array}{cccc}\G(\R,\phi_{I_1^{i_1},2\varepsilon}) & ? & ? & ? \\0 & \G(\R,\phi_{I_2^{i_2},2\varepsilon}) & ? & ? \\\vdots & \vdots & \ddots & ? \\0 & 0 & \dots & \G(\R,\phi_{I_r^{i_r},2\varepsilon})\end{array}\right)  \\ &=
  \left(\begin{array}{cccc}1 & ? & ? & ? \\0 & 1 & ? & ? \\\vdots & \vdots & \ddots & ? \\0 & 0 & \dots & 1 \end{array}\right).
\end{align*}
\normalsize 

Each entry in those matrices is uniquely characterized by one scalar. Hence, we can consider their rank. The left hand side has rank at most equal to the minimum of $r$ and $s$, in particular it is less or equal to $|M(A)|$. The right-hand side has rank $r = |A|$. Therefore we obtain the inequality we wanted. 
\end{pf}

\subsection{Isometry theorem}

In this section, we put together the results proved before to prove that the convolution distance between two sheaves is exactly the same as the bottleneck distance between their graded-barcodes.

\begin{thm}[Isometry]\label{T:DerivedIsometry}
Let $F,G$ be two objects of $\D^b_{\R c}(\kk_\R)$, then : $$ d_C(F,G) =  d_B(\B(F),\B(G)).$$

\end{thm}

\begin{pf}
By Lemma \ref{L:ConvsmallerBottleneck}, there only remains to prove that $ d_C(F,G) \geq d_B(\B(F),\B(G))$, or equivalently, that any $\varepsilon$-interleaving between $F$ and $G$ induces an $\varepsilon$-matching between $\B(F)$ and $\B(G)$.

According to sections 5.3 and 5.4, this interleaving induces a $\varepsilon$-matching between the central, left and right parts of $F$ and $G$, which proves the theorem. 

\end{pf}

\section{Applications}
In this section, we expose some corollaries of the isometry theorem. We start with some explicit computations on an example, showing the fundamentally derived nature of our graded-bottleneck distance.  
Then, we prove that $d_C$ is closed, that is, two sheaves are $\varepsilon$-close if 
and only if they are $\varepsilon$-interleaved, which in particular implies that $d_C$ induces a metric on the isomorphism classes
of $\D^b_{\R c}(\kk_\R)$. We then provide a counter-example of two non constructible sheaves being at convolution distance zero, but which are not isomorphic. These results answer an open question asked by Kashiwara-Schapira in \cite{Kash18} in the one dimensional case. The fact that $d_C$ is closed allows us to consider the set of isomorphism classes of $\D^b_{\R c}(\kk_\R)$ as a  topological metric space. We prove that it is locally path-connected and give a characterization of its connected components.

\subsection{Example : projection from the circle } \label{S:exemplecircle}
We aim here to explain and compute an explicit example that was pointed to us by Justin Curry. It consists of two simple maps from the 
euclidean circle to the real line. 
Understanding this example has been at the origin of our work. 
It is simple yet general enough to exhibit the phenomenons and issues that can happen with the matchings of graded barcodes.

Let $\mathbb{S}^1 = \{(x,y) \in \R^2 \mid x^2 + y^2 = 1 \}$ be the one dimensional circle seen as a sub-manifold in $\R^2$. 
Let  $f : \mathbb{S}^1 \to \R$ be the first coordinate projection and $g : \mathbb{S}^1 \to \R $ 
be the constant map with value zero. Let $F=\Rr f_*\kk_{\mathbb{S}^1}$  and $G=\Rr g_*\kk_{\mathbb{S}^1}$.
Since $\|f - g\|= 1$, the stability theorem by Kashiwara and Schapira~\cite[theorem 2.7]{Kash18} implies :  
$$d_C(F,G) \leq 1  .$$
The CLR decomposition (Definition~\ref{D:CLRdecomposition}) of this two  complexes of sheaves is easy to compute (and depicted 
in the figure below). 
\begin{prop} The complexes $F$ and $G$ have non-zero cohomology spaces at most in degree 0 and 1. Moreover : 
\begin{enumerate}
\item $\Ho^0(F) \simeq \kk_{(-1,1)} \oplus \kk_{[-1,1]}$ and $\Ho^1(F)\simeq0$
\item $\Ho^0(G) \simeq \kk_{\{0\}}$ and  $\Ho^1(G) \simeq \kk_{\{0\}}$
\end{enumerate}

\end{prop}

\begin{center}
\includegraphics[scale=0.8]{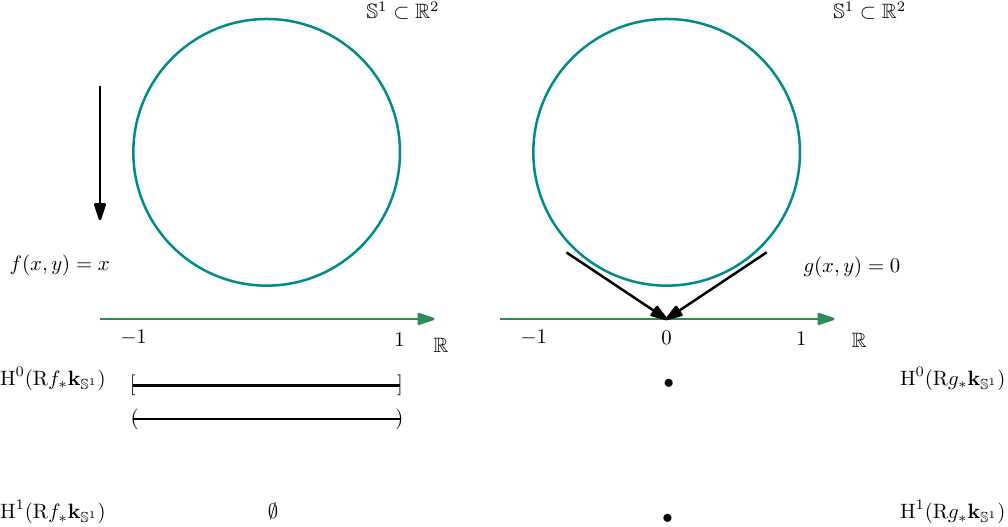}
\end{center}

Hence, $F$ and $G$ are central sheaves and $\B(F) = \{[-1,1]^0,(-1,1)^0\}$, $\B(G)=\{\{0\}^0, \{0\}^1\}$. 
Even in this simple example, there could be \emph{no}
$\varepsilon$-matching between the graded-barcodes if one was working in the ordinary graded category of sheaves. Indeed, $d_C(\kk_{\{0\}}[-1], 0)= + \infty$. However, using our derived notion of interleavings and matching distance we get the expected answer and in fact prove that in this case 
the bound given by the $L_\infty$-norm between the function is optimal. 

 Indeed, let $\sigma$ be the matching defined by :
$$\sigma([-1,1]^0) = \{0\}^0 ~~\text{and}~~\sigma((-1,1)^0) = \{0\}^1.$$

Then we claim that $\sigma$ is a $1$-matching between $\B(F)$ and $\B(G)$. 
Moreover, since the convolution distances between any pair of graded intervals is at least $1$, 
there can not exist an $\varepsilon$-matching between $\B(F)$ and $\B(G)$ for $0\leq \varepsilon < 1$.
Hence we have $F\sim_1 G$ and further
\begin{prop} The convolution distance between $F=\Rr f_*\kk_{\mathbb{S}^1}$  and $G=\Rr g_*\kk_{\mathbb{S}^1}$ is 
$$d_C(\Rr f_*\kk_{\mathbb{S}^1}, \Rr g_*\kk_{\mathbb{S}^1}) = 1.$$
\end{prop}

\subsection{About the closedness of $d_C$ }
In this section we apply our isometry Theorem~\ref{T:DerivedIsometry} to answer an open question of Kashiwara-Schapira 
on the closedness of the convolution distance (see Remark 2.3 of \cite{Kash18}) in the one dimensional case. More precisely, we show that the convolution distance is closed between constructible sheaves over $\R$. We also provide a counter-example to this statement without constructibility assumption. 
\begin{thm}\label{T:closedenessofConvolution}
The convolution distance is closed on $ \D^b_{\R c}(\kk_\R)$. That is, for $F,G \in \D^b_{\R c}(\kk_\R)$ and $\varepsilon \geq 0$ :
$d_C(F,G) \leq \varepsilon \iff F \sim_\varepsilon G. $

\end{thm}

We start with the following easy  lemma, whose proof is left to the reader.

\begin{lem}
Let $I^i,J^j$ two graded intervals (possibly empty, we set $\kk_\emptyset = 0$) and $\varepsilon \geq 0$. Then : $$d_C(\kk_I[-i],\kk_J[-j]) \leq \varepsilon \iff \kk_I[-i] \sim_\varepsilon \kk_J[-j]  $$
\end{lem}

\begin{pf}[of the theorem]

Suppose $d_C(F,G)\leq \varepsilon$. Then by definition there exists a decreasing sequence $(\varepsilon_n)$ such that $\varepsilon_n \rightarrow \varepsilon $ when $n$ goes to infinity and for every $n\in\N$, $F\sim_{\varepsilon_n} G$. For simplicity of the proof, we will assume the graded-barcodes of $F$ and $G$ to be finite, but the proof generalizes to the locally finite case. Then by applying the isometry theorem, for $n\geq 0$, there exists a $\varepsilon_n$ matching $\sigma_n : \B(F) \rightarrow \B(G)$. 

Now by finiteness of the graded-barcodes, the set of matchings between $ \B(F)$ and $ \B(G)$ is finite. Hence, we can extract from $(\sigma_n)$ a constant sequence, say $(\sigma_{\varphi(n)})$. Applying lemma 4.1 and making $n$ going to infinity, we see that $\sigma := \sigma_{\varphi(0)}$ is an $\varepsilon$-matching between $ \B(F)$ and $ \B(G)$.

\end{pf}

\begin{remark}
One must observe that in the case of persistence modules, the interleaving distance is \emph{not} closed. There exists some ephemeral modules at distance 0 from 0 : consider the one parameter persistence module $\kk^{\{0\}}$ (keeping notations of section 5.2). To avoid this issue, Chazal,
Crawley-Boevey and de Silva introduced the observable category of persistence modules $\text{Obs}(\Pe(\kk_\R))$ in \cite{Chaz16}. It is defined as the quotient category of $\Pe(\kk_\R)$ by the full sub-category of ephemeral persistent modules, which has objects $M\in \Pe(\R)$ such that $M(s<t)=0$ for every $s<t \in \R$. By construction, the interleaving distance on $\Pe(\kk_\R)$ induces a closed metric on $\text{Obs}(\Pe(\kk_\R))$. Note that this construction has since been generalized by the first author in \cite{BerPetit}.

\end{remark}

\begin{cor}
The functors $\Psi^j_R : \D^j_R \to \Pe(\R)$ (see proposition \ref{P:ExistenceofPsi}) induces an isometric equivalence of categories between $\D^j_R$ and $\text{Obs}(\Pe(\R))^{op}$.
\end{cor}

We now explicit a counter-example to the closedness of $d_C$ without constructibility assumptions. More precisely, we will construct two sheaves $F,G \in \D^b(\kk_\R)$ such that $d_C(F,G)= 0$ but $F \not \simeq G$. We consider the sets $X = \Q \cap [0,1]$ and $Y= \sqrt{2}\Q \cap [0,1] = \{\sqrt{2} q \mid q\in \Q\}\cap [0,1]$. 

\begin{prop}\label{prop:approxuniforme}
There exists a sequence of functions $(r_n)_{n\in\Z_{> 0}}$ from $X$ to $Y$ satisfying: \begin{enumerate}
    \item for any $n \in \Z_{> 0}$, $r_n : X \to Y$ is bijective,
    \item $\sup_{x \in X} |r_n(x) - x| \xrightarrow[n \to + \infty]{} 0$.
\end{enumerate}
\end{prop}

\begin{pf}
 Let $n\in \Z_{> 0}$. We define $r_n$ piecewise on $[0,1-\frac{1}{n}]\cap \Q$ and $]1-\frac{1}{n},1]\cap \Q$. For $q\in [0,1-\frac{1}{n}]\cap \Q$, we set $r_n(q) = \frac{\sqrt{2}}{\lceil \sqrt{2} \rceil_{10}(n)}q$, with $\lceil \sqrt{2} \rceil_{10}(n) = \frac{\lceil10^n \sqrt{2}\rceil }{\lceil{10^n}\rceil}$ the $n$-th ceil decimal approximation of $\sqrt{2}$. Then ${r_n}_{|[0,1-\frac{1}{n}]\cap \Q}$ is injective, and 
 
 \begin{align*}
r_n\left ([0,1-\frac{1}{n}]\cap \Q\right ) &= \left \{\sqrt{2}q \mid q \in \left [0, \frac{1-1/n}{\lceil \sqrt{2}  \rceil_{10}(n)}\right ]\cap \Q  \right\} \\
& = Y\cap \left[0,\sqrt{2}\frac{1-1/n}{\lceil \sqrt{2}  \rceil_{10}(n)} \right] \\ 
 &\subsetneq Y. 
 \end{align*} 
 
 Now, since $]1-\frac{1}{n},1]\cap \Q$ and $Y \big \backslash r_n\left ([0,1-\frac{1}{n}]\cap \Q\right )$ are both infinite subsets of $\Q$, there exists a bijection 
 
 $$\varphi_n :~ \left]1-\frac{1}{n},1 \right ]\cap \Q \xrightarrow[]{\sim} Y \big \backslash r_n\left ([0,1-\frac{1}{n}]\cap \Q\right ). $$
 
 We define ${r_n}_{| ]1-\frac{1}{n},1]\cap \Q } = \varphi_n$. Then ${r_n}_{| ]1-\frac{1}{n},1]\cap \Q }$ is injective and $$r_n(]1-\frac{1}{n},1]\cap \Q) = Y \big \backslash r_n\left ([0,1-\frac{1}{n}]\cap \Q\right ).$$
 
 Finally, $r_n$ is indeed a bijective function from $X$ to $Y=\sqrt{2}\Q \cap [0,1]$. 
 
\center{\includegraphics[scale = 4.9]{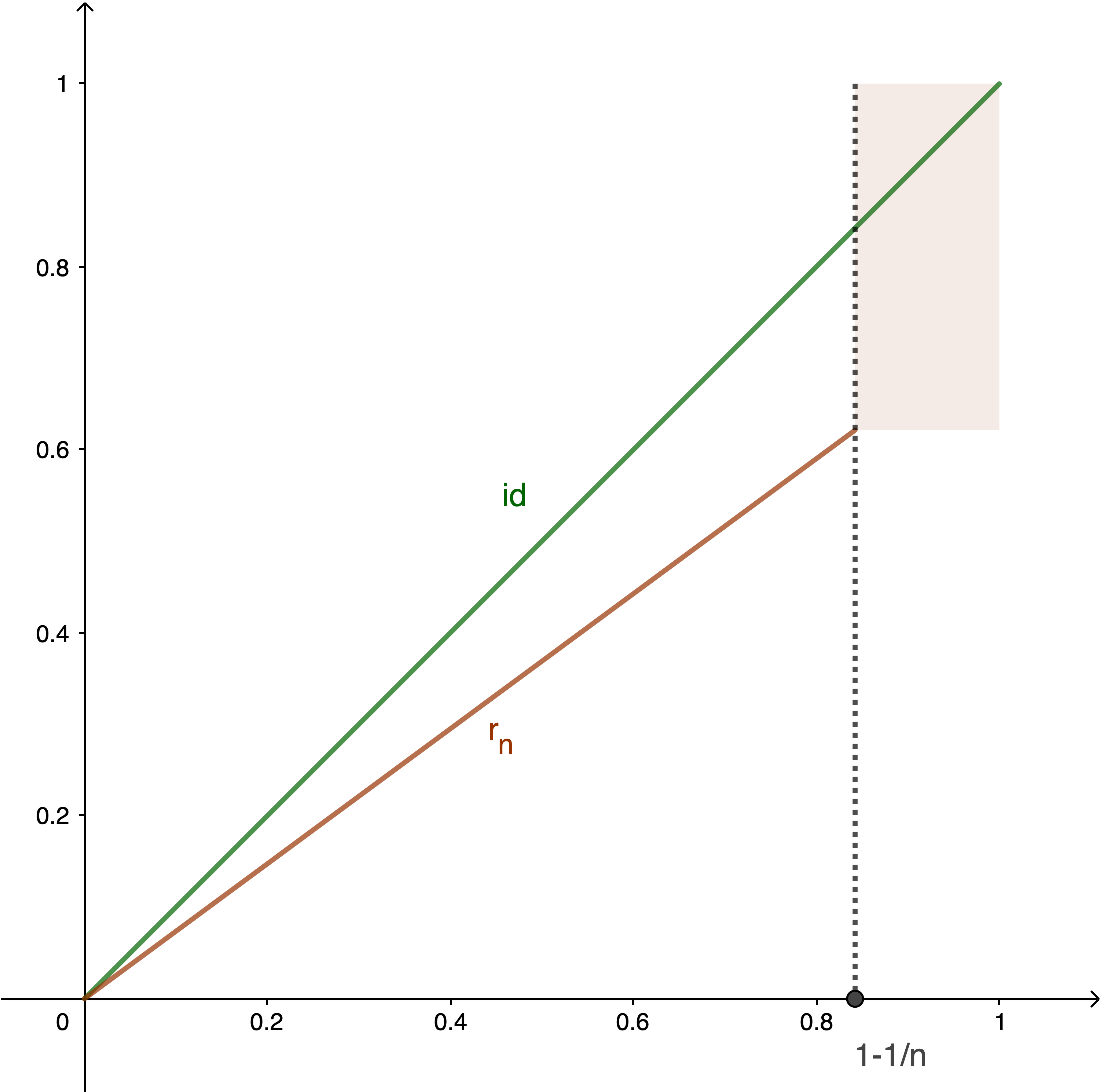}}

\emph{Graphical representation of} $r_n$
\flushleft
 Observe that : 
 
 $$\sup_{x\in X} |r_n(x) - x| = \max \left (\sup_{x\in X\cap [0,1-\frac{1}{n}]} |r_n(x) - x|, \sup_{x\in X\cap [1-\frac{1}{n},1]} |r_n(x) - x|  \right ).$$
 
 The first term of the maximum is worth $\left ( \frac{\sqrt{2}}{\lceil \sqrt{2} \rceil_{10}(n)} - 1 \right )\left (1 - 1/n \right )$, and the second term is bounded from above by the diameter of the interval $[ \sqrt{2}\frac{1-1/n}{\lceil \sqrt{2}  \rceil_{10}(n)} , 1]$ which is worth $1- \sqrt{2}\frac{1-1/n}{\lceil \sqrt{2}  \rceil_{10}(n)} $. Since both of these terms go to $0$ as $n$ goes to infinity, we deduce the desired property :  $$ \sup_{x \in X} |r_n(x) - x| \xrightarrow[n \to + \infty]{} 0.$$
 
 \end{pf}

\begin{prop}\label{prop:suminterleaving} Let $(F_i)_{i\in I}$ and $(G_j)_{j\in J}$ two families of objects of $\D^b(\kk_\V)$. Assume that there exists a bijective function $\sigma: I \to J $, and $\varepsilon\geq 0$ such that for all $i\in I$, $d_C(F_i,G_{\sigma(i)})\leq \varepsilon$. Then : $$d_C \left (\bigoplus_{i\in I} F_i, \bigoplus_{j\in J} G_j  \right) \leq \varepsilon.  $$ 

\end{prop}

\begin{pf}
Let $\varepsilon' > \varepsilon$ and $i\in I$. Then by assumptions, there exists $\varepsilon'$-interleaving morphisms between $F_i$ and $G_{\sigma(i)}$, $\varphi_i : F_i \star K_{\varepsilon'} \to G_{\sigma(i)}$ and $\psi_i : G_{\sigma(i)}\star K_{\varepsilon'} \to F_i$.
Since $-\star K_{\varepsilon'}$ is a left-adjoint functor, it commutes with arbitrary colimits. Therefore, by taking direct sums of the previous $\varepsilon'$-interleaving morphisms,  we get $\varepsilon'$-interleaving morphisms between $\bigoplus_{i\in I} F_i$ and $\bigoplus_{i \in I} G_{\sigma(i)} \simeq \bigoplus_{j \in J} G_{j}$, which proves the result.
\end{pf}

Let $F = \bigoplus_{x\in X} \kk_{\{ x\}}$ and $G = \bigoplus_{y\in Y} \kk_{\{ y\}}$.

\begin{prop}
$F$ is not isomorphic to $G$ and  $d_C(F,G)=0$.
\end{prop}

\begin{pf} 
$F$ and $G$ cannot be isomorphic since $F_1 \simeq \kk$ and $G_1 \simeq 0$.

 Let $r_n : X \to Y$ be as in proposition \ref{prop:approxuniforme}. Using proposition \ref{prop:suminterleaving}, and the fact that, for $x,y\in \R$, $d_C(\kk_{\{x\}},\kk_{\{y\}}) = |x-y|$, we obtain that, for any $n \in \Z_{>0}$:
 
 $$d_C(F,G) \leq \sup_{x\in X} |r_n(x)-x|.$$
 
Taking the limit as $n$ goes to infinity, we deduce that $d_C(F,G)=0$.
\end{pf}

\subsection{Description of the connected components of $\D^b_{\R c}(\kk_\R)$}
In this section, we study the connected components of $\D^b_{\R c}(\kk_\R)$ with respect to the metric. In order to make sense of this, we introduce the \emph{small} category $Barcode$ which is a combinatorial description of $\D^b_{\R c}(\kk_\R)$ and is equivalent to  $\D^b_{\R c}(\kk_\R)$. The category $Barcode$ is shown to be skeletal (any two isomorphic objects are equal), and is equipped with the graded bottleneck distance. From $Barcode$ we thus obtain an extended metric space, which will be proven to be locally path-connected, see Theorem~\ref{T:pathcomponent}. To do so, we first prove an interpolation lemma in the same fashion as Chazal et al. \cite[Theorem 3.5]{ChazSilv16}, which stands that if two sheaves are $\varepsilon$-interleaved, there exists a $1$-lipschitz path in $D^b_{\R c}(\kk_\R)$ between them.

\begin{lem}[Interpolation]\label{L:interpolation}
Let $F,G \in \D^b_{\R c}(\kk_\R)$ be such that $F\sim_{\varepsilon} G$ for some $\varepsilon \geq 0$. Then there exists a family of sheaves $(U_t)_{t\in[0,\varepsilon]}$ in $\D^b_{\R c}(\kk_\R)$ such that :

\begin{enumerate}
\item $U_0 = F$ and $U_\varepsilon = G$.
\item For $t\in [0,\varepsilon]$, $d_c(F,U_t) \leq t$ and $d_C(G,U_t) \leq \varepsilon - t$.
\item For $(t,t')\in [0,\varepsilon]^2$, $d_C(U_t,U_{t'}) \leq |t-t'|$.
\end{enumerate}

\end{lem}

\begin{pf} Let $F\star \K_\varepsilon \stackrel{\varphi}{\longrightarrow} G$ and $G\star \K_\varepsilon \stackrel{\psi}{\longrightarrow} F$ be the interleaving morphisms between $F$ and $G$.

We start by constructing $U_t$ for $t\in [0,\frac{\varepsilon}{2}]$. 
The interleaving morphism and the canonical maps in $\D^b_{\R c}(\kk_\R)$  give (by Proposition~\ref{P:propertiesofconvolution} 
and~\eqref{eq:propertiesofconvolution})  the following diagram $\mathbb{D}_t$:

$${\xymatrix{
    G \star K_{t - \varepsilon} \ar[dd]_{\phi_{G,2t}\star K_{t-\varepsilon }} \ar[rrdd]^(.3){\;\;\psi\star K_{t-\varepsilon}\qquad\qquad} & & F \star K_{-t} \ar[dd]^{\phi_{F,2\varepsilon-2t}\star K_{-t}} \ar[ddll]^(.3){\varphi \star K_{-t}\qquad} \\
    & & \\
    G \star K_{-t-\varepsilon} & & F\star K_{t-2\varepsilon}
  }.}$$ 
  Taking resolutions in $\Mod(\kk_\R)$, one can assume this diagram is actually given by a  diagram still denoted 
  $\mathbb{D}_t$ in $C(\Mod(\kk_\R))$ which we assume from now on. 
 One can note that this diagram defines  two maps $\theta_t$, $ \tilde{\phi}_t:(G \star K_{t - \varepsilon}) \oplus (F \star K_{-t})  \longrightarrow 
(G \star K_{-t-\varepsilon}) \oplus (F\star K_{t-2\varepsilon})$
given by $$(x,y) \stackrel{\theta_t}\longmapsto (\varphi \star K_{-t}(y),\psi\star K_{t-\varepsilon}(x))\mbox{ and } (x,y)\stackrel{\tilde{\phi}_t}\longmapsto (\phi_{G,2t}\star K_{t-\varepsilon }(x), 
\phi_{F,2\varepsilon-2t}\star K_{-t}(y)).$$ 
The limit $\varprojlim \mathbb{D}_t$ of the diagram is precisely (isomorphic  to)  
the equalizer of the two maps and thus to the kernel $\ker(\theta_{t} - \tilde{\phi}_{t})$ of their difference. It is now enough to define  ${U}_t := \mathrm{ho}\varprojlim \mathbb{D}_t$ to be the homotopy limit in (the model category of sheaves~\cite{Crans95}) $C(\Mod(\kk_\R))$ of the diagram $\mathbb{D}_t$ together with the canonical maps from the kernel to $U_t$ and from $U_t$ to $F\star K_{-t}$ given by the diagram to conclude. This is what we do below using an explicit model $\tilde{U}_t$ without further reference to or use of  homotopy limit. 

%Since we are dealing with a diagram in $C(\Mod(\kk_\R))$ that we wish to see in the derived category we essentially only need to replace the limit by its homotopy limit.  
%Namely, we define $\tilde{U}_t := \mathrm{ho}\varprojlim \mathbb{D}_t$ to be the homotopy limit in (the model category of sheaves~\cite{Crans95}) $C(\Mod(\kk_\R))$ of the diagram $\mathbb{D}_t$.  For the reader who wish to avoid the use of the technique of homotopy limits in model category (see~\cite{DuggerHocolim}),
Let us denote $A_t:= (G \star K_{t - \varepsilon}) \oplus (F \star K_{-t})$ and $B_t:=(G \star K_{-t-\varepsilon}) \oplus (F\star K_{t-2\varepsilon})$. We define 
$\tilde{U}_t:= \mathrm{cocone}(A_{t} \stackrel{}\longrightarrow B_{t})$ that is the complex of sheaf $A_{t} \oplus B_{t}[-1]$ endowed with the differential $\left(\begin{array}{cc} d^A & 0 \\ \theta_{t} - \tilde{\phi}_{t}& -d^B \end{array}\right)$. 

We need to prove that $\tilde{U}_t$ is $t$-interleaved with $F$.
Note that the canonical projection $A_{t} \oplus B_{t}[-1] \to A_{t}$ gives a chain map $\tilde{U}_{t} \to A_{t}$. Since $A_{t}= (G \star K_{t - \varepsilon}) \oplus (F \star K_{-t})$ we can compose the latter with the projection on either factors of $A_{t}$ as well, and in particular we have 
$\tilde{f}:\tilde{U}_t\to A_{t}\to F \star K_{-t}$ and hence  (by proposition~\ref{P:propertiesofconvolution}) the map
\begin{equation}
 \tilde{U}_t\star K_{t} \stackrel{f}\longrightarrow F. 
\end{equation}

We now need to define a map $g: F\star K_t \to \tilde{U}_t$. First note that the first summand inclusion of $A_{t}$ into  $\tilde{U}_{t}$ is not a chain map  but the composition 
$\iota: \ker(\theta_{t} - \tilde{\phi}_{t}) \to A_{t}\hookrightarrow \tilde{U}_t$ is a morphism  in $C(\Mod(\kk_\R))$. Now,  the interleaving map $\varphi: F\star K_{\varepsilon}\to G$ induces the map 
\begin{equation}\label{eq:secpondpartinterleaving} F\star K_t \stackrel{(\varphi\star K_{t-\epsilon}, \phi_{F,2t}\star K_t)}\longrightarrow 
(G \star K_{t - \varepsilon}) \oplus (F \star K_{-t})\end{equation} which makes the following diagram 
$${\xymatrix{ & F\star K_t  \ar[rd]^{\phi_{F,2t}\star K_t} \ar[ld]_{\varphi\star K_{t-\epsilon}}&\\
    G \star K_{t - \varepsilon} \ar[dd]_{\phi_{G,2t}\star K_{t-\varepsilon }} \ar[rrdd]^(.3){\;\;\psi\star K_{t-\varepsilon}\qquad\qquad} & & F \star K_{-t} \ar[dd]^{\phi_{F,2\varepsilon-2t}\star K_{-t}} \ar[ddll]^(.3){\varphi \star K_{-t}\qquad} \\
    & & \\
    G \star K_{-t-\varepsilon} & & F\star K_{t-2\varepsilon}
  }}$$ commutative since $\varphi$, $\psi$ defines a $\varepsilon$-interleaving. This implies that the map~\eqref{eq:secpondpartinterleaving} factors through  $\varprojlim \mathbb{D}_t\cong  \ker(\theta_{t} - \tilde{\phi}_{t}) $ and hence we get 
  the  map \[g: F\star K_t\longrightarrow \varprojlim \mathbb{D}_t\stackrel{\iota}\longrightarrow  \tilde{U}_t\] in $\D^b_{\R c}(\kk_\R)$. The maps $f$ and $g$ gives us the required interleaving because $\varphi$ and $\psi$ are.

For $t\in ]\frac{\varepsilon}{2}, \varepsilon]$, we construct $U_t$ in a similar fashion by intertwining the roles of $F$ and $G$ in the diagram $C(\Mod(\kk_\R))$.

Let $\Delta_\varepsilon = \{(x,y) \in \R^2 \mid 0 \leq {y-x} \leq \varepsilon \}$ be equipped with the standard product order of $\R^2$ : $(x,y)\leq (x',y') \iff x\leq x'$ and $y\leq y' $. Observe that the mapping : $$\Delta_\varepsilon \ni (x,y) \rightsquigarrow U_{y-x}\star K_{-x-y}$$
induces a well defined functor $(\Delta_\varepsilon, \leq) \longrightarrow \D^b_{\R c}(\kk_\R)$ whose restriction to the poset $\{(x,y) \in \R^2 \mid y-x = t\}$ is the functor : $(x,y) \longrightarrow U_t\star K_{-x-y}$ with internal maps given by the natural morphisms $(\phi_{U_t,\varepsilon})$. Hence, for $ \varepsilon \geq t,t'\geq 0$, $U_t$ and $U_{t'}$ are $|t-t'|$ interleaved.
\end{pf}

We will now define the category $Barcode$ we mentioned earlier.

We first setup notations and terminology for (graded) intervals (with multiplicity).
Let $\text{Int}(\R)$ be the set of intervals of $\R$ and $p_1$,  $p_2$ be the two first coordinate projections of  $\text{Int}(\R) \times \Z \times \Z_{\geq 0}$. Let $\mathbb{B}$ be a subset of $ \text{Int}(\R) \times \Z \times \Z_{\geq 0}$. Then  $\mathbb{B}$ is said to be 
\begin{itemize}\item \emph{locally finite} if $p_1(\mathbb{B})\cap K$ is finite for all compact subsets of $\R$; \item \emph{bounded} if $p_2(\mathbb{B}) \subset \Z$ is bounded;
\item \emph{well-defined} if the fibers of the projection $(p_1,p_2)$ have cardinality at most $1$.
\end{itemize}
In a triple $(I,j,n) \in\mathbb{B}$, the first integer will stand for the degree on which the interval $I$ is seen and the second non-negative integer $n$ stands for its multiplicity.

\begin{defi}\label{D:CatBarcode}The category $Barcode$ has objects the set
$$\Obj({Barcode}) = \{\mathbb{B} \subset \text{Int}(\R) \times \Z\times \Z_{\geq 0} \mid \mathbb{B} ~\text{is bounded, locally finite and well-defined} \}. $$

\noindent For any $\mathbb{B}$ and $\mathbb{B}' \in {Barcode}$,  the set of their morphisms is
$$ \Hom_{{Barcode}}(\mathbb{B},\mathbb{B}') =  \prod_{\substack{(I,j,n)\in \mathbb{B} \\ (I',j',n') \in \mathbb{B}'}}\Hom_{\D^b_{\R c}(\kk_\R)} \left (\kk_I^n[-j] , \kk_{I'}^{n'}[-j'] \right ).
$$
\end{defi}

We define the composition in {Barcode} so that the mapping : $$\iota : \Obj({Barcode})\ni \mathbb{B}  \mapsto \bigoplus_{(I,j,n)\in \mathbb{B}} \kk_I^n[-j] \in \Obj(\D^b_{\R c}(\kk_\R))  $$ becomes a fully faithful functor : 
$$\iota :  {Barcode} \longrightarrow \D^b_{\R c}(\kk_\R).$$

Note that this is possible only because the objects of $Barcode$ are locally finite. Theorems \ref{T:KSstructure} and \ref{T:KSdecomposition} assert that $\iota$ is essentially surjective, therefore is an equivalence. We also deduce from these theorems that $Barcode$ is a skeletal category: it satisfies for any $\mathbb{B},\mathbb{B}' \in {Barcode}$, $$\mathbb{B} \simeq \mathbb{B}' \text{~if and only if~} \mathbb{B} = \mathbb{B}'. $$ 

The notion of equality is well-defined here since $\Obj({Barcode})$ is a set. Therefore $\iota$ identifies its image as a skeleton of $\D^b_{\R c}(\kk_\R)$, a full-subcategory which is dense and skeletal.

Moreover, Theorems \ref{T:KSstructure} and \ref{T:KSdecomposition} allows us to  equip the set $\Obj({Barcode})$ with the graded-bottleneck distance (definition \ref{D:bottleneckgeneral}). The derived isometry theorem \ref{T:DerivedIsometry} implies that, for any $\mathbb{B},\mathbb{B}' \in {Barcode}$, one has : 

$$d_C(\iota(\mathbb{B}),\iota(\mathbb{B}')) = d_B(\mathbb{B},\mathbb{B}'). $$

\begin{thm}\label{T:pathcomponent}
The following assertions hold: 

\begin{enumerate}
    \item $(\Obj({Barcode}),d_B)$ is an extended metric space,
    \item $(\Obj({Barcode}),d_B)$ is locally path-connected.
\end{enumerate}
\end{thm}

\begin{pf}
\begin{enumerate}
    \item The fact that $d_B$ is a pseudo-extended metric is inherited from the properties of $d_C$ (proposition \ref{P:propertiesconvdistance}) by the derived isometry theorem. Moreover, if $d_B(\mathbb{B},\mathbb{B}')=0$ then $\mathbb{B}=\mathbb{B}'$ by theorem \ref{T:closedenessofConvolution}.
    
    \item We will prove that open balls are path-connected, that is, any two barcodes at finite distance can be connected by a continuous path. Let $\mathbb{B}_0$ and $\mathbb{B}_\varepsilon$ in ${Barcode}$ such that $d_B(\mathbb{B}_0,\mathbb{B}_\varepsilon) = \varepsilon$. According to the interpolation lemma \ref{L:interpolation}, there exists a family of objects $(F_t)_{t\in [0,\varepsilon]}$ of $\D^b_{\R c}(\kk_\R)$ such that $F_0 = \iota(\mathbb{B}_0)$, $F_\varepsilon = \iota(\mathbb{B}_\varepsilon)$, and for any $t,t'\in [0,\varepsilon]$, $d_C(F_t,F_{t'})\leq |t-t'|$. Given $t\in [0,\varepsilon]$, define $\mathbb{B}_t$ to be the graded-barcode of $F_t$. Then, it is clear thanks to the derived isometry theorem that $(t \mapsto \mathbb{B}_t)$ defines a 1-lipschitz path between $\mathbb{B}_0$ and $\mathbb{B}_\varepsilon$.
    
\end{enumerate}

\end{pf}

\subsection{Algorithmic remarks on computing one best matching}

The formulation of the convolution distance as a matching distance we obtained in Section~5 turns the computation of an algebraic problem into minimizing the cost of a matching, which is of combinatorial nature. This is in fact a variant of a very classical problem of linear programming, for which there exists an abundant literature that can be solved in polynomial time using the \emph{Hungarian algorithm} \cite{Kuhn09}. Hence, distances in $\D^b_{\R c}(\kk_\R)$ \emph{can be implemented} in a computer and \emph{computed}.

\bibliographystyle{alpha}
\newcommand{\etalchar}[1]{$^{#1}$}

\end{document}